\documentstyle{amsppt}
\frenchspacing
\catcode`\@=11
\def\Imm{\operatorname{Imm}}
\def\det{\operatorname{det}}
\def\id{\operatorname{id}}
\def\dim{\operatorname{dim}}
\def\R{\Bbb R}
\def\quad{\null\hskip1pc\relax}
\def\qquad{\null\hskip2pc\relax}

\hsize=38pc
\pageheight{20cm}
\newdimen\mathindent	\mathindent=3pc
\newdimen\disp@box@width \disp@box@width=35pc 
\newdimen\disp@width	\disp@width=32pc	
\def\eqno#1{\null\hfill\llap{{\rm#1}\unskip}}
\def\@fl#1$${\def\paritem[##1]{%
\llap{\hbox to 4pc{\hbox to \paritemwd{\hfil\rm ##1}\hfil}}} 
\hbox to\disp@width{$\vcenter{\parindent=0pc\hsize\disp@box@width
\let\\=\eqn@cr$\displaystyle#1$\hfil}$\hss}$$
\ignorespaces} \everydisplay{\@fl}
\def\eqn@cr{\relax$\par\vskip5pt$\displaystyle} \def\tag#1{\eqno{#1}}
\def\spreadmatrixlines#1{\spreadmlines@#1\relax} 

\def\aligned{\vcenter\aligned@}

\def\gathered{\null\,\vcenter\bgroup\vspace@\Let@ 
\ifinany@\else\openup\jot\fi\ialign
\bgroup\strut@$\m@th\displaystyle{##}$\hfil
&&\quad\strut@$\m@th\displaystyle{##}$\hfil\crcr} 
\def\endgathered{\crcr\egroup\egroup}

\let\le=\leqslant	\let\ge=\geqslant
\let\leq=\leqslant	
	
\newdimen\paritemwd
\paritemwd=2pc
\def\paritem[#1]{\par\noindent\hbox to\paritemwd{\hss \rm #1}}
\catcode`\@=\active

\topmatter

\title
INVARIANTS OF VELOCITIES, AND HIGHER ORDER GRASSMANN BUNDLES 
\endtitle

\author
DAN RADU GRIGORE AND DEMETER KRUPKA
\endauthor
\footnote" "{Research of the first author supported by Grant No. 871/95 of
the Ministry of Education and Youth of 
the Czech Republic. Research of the second author supported by Grant No.
201/96/0845 of the Czech Grant Agency.} 

13, 74601 Opava, Czech Republic 


\endtopmatter

\document

invariants of $(r,n)$-velocities 

\hrule
\vskip 1mm

\noindent{\bf Abstract}

\vskip 3mm
An
$(r,n)$-velocity
is an
$r$-jet
with source at
$0 \in \R^n$,
and target in a manifold
$Y$. An
$(r,n)$-velocity
is said to be regular, if it has a representative which is an immersion at
$0 \in \R^{n}$.
The manifold
$T^{r}_{n}Y$
of
$(r,n)$-velocities
as well as its open,
$L^{r}_{n}$-invariant,
dense submanifold
$\Imm T^{r}_{n}Y$
of regular
$(r,n)$-velocities,
are endowed with a natural action of the differential group $L^{r}_{n}$
of invertible
$r$-jets
with source and target
$0 \in \R^{n}$.
In this paper, we describe all continuous, $L^{r}_{n}$-invariant,
real-valued functions on
$T^{r}_{n}Y$
and
$\Imm T^{r}_{n}Y$.
We find local bases of
$L^{r}_{n}$-invariants
on
$\Imm T^{r}_{n}Y$
in an explicit, recurrent form. To this purpose, higher order Grassmann
bundles are considered as the 
corresponding quotients $P^{r}_{n}Y = \Imm T^{r}_{n}Y/L^{r}_{n}$, and their
basic properties are studied. 
We show that nontrivial $L^{r}_{n}$-invariants
on
$\Imm T^{r}_{n}Y$
cannot be continuously extended onto
$T^{r}_{n}Y$.
\vskip 3mm

\noindent{\it Subj. Class.:}
Differential invariants; Grassmann bundles



\vskip 3mm
\hrule

\heading 1. Introduction \endheading

By a {\it velocity} one usually means the derivative of a curve in a smooth
manifold
$Y$
at a point, or, which is the same, the tangent vector to this curve at a
point $y\in Y$.
Equivalently, such a velocity is a
$1$-jet
with {\it source} at the origin
$0 \in \R^{n}$
and {\it target} in
$Y$.
Generalizing this concept one may define an $(r,n)$-{\it velocity}
as an
$r$-jet
with source
$0 \in \R^{n}$
and target in
$Y$.

If such an
$r$-jet
can be represented by an immersion of a neighbourhood of the origin $0\in
\R^{n}$
into
$Y$,
it is called {\it regular}, and we speak of a {\it regular} $(r,n)$-{\it
velocity}.

Concepts of this kind, i.e., the
$r$-jets
of differentiable mappings between smooth manifolds, have been introduced
in the fiftieth by Ehresmann 
(see references in [7]), and have become the basic concepts of the theory
of differential invariants, 
and the theory of natural bundles and operators 
(see Kol\'a\v{r}, Michor, Slov\'ak [7], D. Krupka, Jany\v{s}ka [9], D.
Krupka [11], Nijenhuis [14], and the 
references therein).  
It should be pointed out, however, that the problem of finding invariants
of velocities and the corresponding 
problem of describing the structure of the space of higher order velocities
has not been touched in the existing 
monographs on differential invariants and natural bundles [7], [9]. 

The set of
$(r,n)$-velocities
on a smooth manifold is a smooth manifold endowed with a right action of
the differential group
$L^{r}_{n}$
of invertible
$r$-jets
with source and target
$0 \in \R^{n}$.
The purpose of this paper is to characterize all continuous scalar
invariants of this action, i.e. all real-valued 
functions defined on open subsets of $T^{r}_{n}Y$,
which are constants on the
$L^{r}_{n}$-orbits.
Instead of formulating and solving equations for invariant functions we use
a different, more powerful method 
based on considering the quotient space of the open, dense subspace of
$T^{r}_{n}Y$,
formed by regular
$(r,n)$-velocities.
The corresponding orbit space is then called the $(r,n)$-{\it Grassmann
bundle}.
It is a fiber bundle over
$Y$
whose type fiber is the
$(r,n)$-{\it Grassmannian}.
The canonical quotient projection of the manifold of regular
$(r,n)$-velocities
onto the orbit space is the {\it basis of invariants} of
$(r,n)$-velocities.
Geometric structures of this kind as well as their invariants have been
studied by M. Krupka [12], [13].

Thus, to find all
$L^{r}_{n}$-invariants
it is enough to find the projection of the manifolds of regular higher
order velocities onto the higher order 
Grassmann bundle. We note that an analogous method has been applied to the
problem of finding 
$GL_{n}(\R)$-invariants
of a linear connection [10].

Basic concepts of the {\it first order} Grassman bundles has been applied
in mathematical physics, 
and the parametrization independent variational theory (see e.g. Dedecker
[1], Grigore [3], Grigore, Popp [4], 
Horv\'athy [5], Klein [6]). {\it Higher order} Grassmann bundles have
become natural underlying structures for 
the geometric theory of partial differential equations (Krasilschik,
Vinogradov, Lychagin [8]).

\heading 2. Higher order velocities \endheading 

Throughout this paper,
$m$, $n\ge 1$
and
$r\ge 0$
are integers such that
$n\le m$,
and
$Y$
is a smooth manifold of dimension
$n+m$.

By an
$(r,n)$-{\it velocity}
at a point
$y\in Y$
we mean an
$r$-jet
$J^{r}_{0}\zeta$
with {\it source}
$0 \in \R^{n}$
and {\it target}
$y=\zeta (0)$.
The set of
$(r,n)$-velocities
at
$y$
is denoted by
$J^{r}_{(0,y)}(\R^{n},Y)$.
Further, we denote
$$
T^{r}_{n}Y = \bigcup_{y \in Y} J^{r}_{(0,y)}(\R^{n},Y), $$
and define surjective mappings
$\tau^{r,s}_n:T^{r}_{n}Y \to T^{s}_{n}Y$, where
$0\le s \le r$,
by
$
\tau^{r,s}_{n}(J^{r}_{0}\zeta) = J^{s}_{0}\zeta. $
Recall that the set
$T^{r}_{n}Y$
has a smooth structure defined as follows. Let $(V,\psi)$,
$\psi =(y^{A})$,
be a chart on
$Y$..
Then the {\it associated chart}
$(V^{r}_{n},\psi^{r}_{n})$
on
$T^{r}_{n}Y$
is defined by
$
V^{r}_{n} = (\pi^{r,0}_{n})^{-1}(V),\psi^{r}_{n} =
(y^{A},y^{A}_{i_{1}},y^{A}_{i_{1}i_{2}},
\dots, y^{A}_{i_{1}i_{2} \dots i_{r}}),
$
where
$
1\le i_{1} \le i_{2} \le \dots \le i_{r} \le n, $
and for every
$J^{r}_{0} \zeta \in V^{r}_{n}$,

$$
y^{A}_{i_{1}i_{2} \dots i_{k}}(J^{r}_{0}\zeta) = D_{i_{1}}D_{i_{2}} \dots
D_{i_{k}} (y^{A}\zeta)(0), 
\quad 0 \le k\le r.
\tag{(2.1)}
$$
The set
$T^{r}_{n}Y$
endowed with the smooth structure defined by the associated charts is
called the {\it manifold of}
$(r,n)$-{\it velocities}
over
$Y$.

The equations of the mapping
$\tau^{r,s}_{n}: T^{r}_{n}Y \to T^{s}_{n}Y$ in terms of the associated
charts are given by $
y^{A}_{i_{1}i_{2}\dots i_{k}}\circ\tau^{r,s}_{n}(J^{r}_{0}\zeta) =
y^{A}_{i_{1}i_{2} \dots i_{k}}(J^{r}_{0}\zeta), $
where
$
0\le k\le s.
$
In particular, these mappings are all submersions. 

Let
$\text{tr}_{t}$
denote the translation
$t'\to t'+t$ of $\R^{n}$.
If
$\gamma$
is a smooth mapping of an open set
$U \subset \R^{n}$
into
$Y$,
then for any
$t\in U$,
the mapping
$t' \to \gamma \circ\text{tr}_t(t')$
is defined on a neighbourhood of the origin $0 \in \R^{n}$
so that the
$r$-jet
$J^{r}_{0}(\gamma\circ\text{tr}_{t})$
is defined. The mapping

$$
U \ni t \to (T^{r}_{n} \gamma)(t) = J^{r}_{0}(\gamma\circ\text{tr}_{t}) \in
T^{r}_{n}Y
\tag{(2.2)}
$$
is called the
$r$-{\it prolongation},
or simply the prolongation of
$\gamma$.
Since
$
y^{A}_{i_{1}i_{2} \dots i_{k}} \circ T^{r}_{n}\gamma(t) =
D_{i_{1}}D_{i_{2}} \dots D_{i_{k}}(y^{A}(\gamma 
\circ\text{tr}_{t}))(0) $
and
$
D_{i_{1}}(y^{A}(\gamma \circ\text{tr}_{t}))(t') = D_{i_{1}}(y^{A}\gamma
)(t'+t) $,
we get for its chart expression

$$
y^{A}_{i_{1}i_{2} \dots i_{k}} \circ T^{r}_{n}\gamma (t) =
D_{i_{1}}D_{i_{2}}\dots D_{i_{k}}(y^{A}\gamma )(t). 
\tag{(2.3)}
$$

Assume that we have an element
$J^{r}_{0}\zeta \in T^{r}_{n}Y$.
$J^{r}_{0}\zeta$
defines the tangent mapping
$T_{0}T^{r-1}_{n}\zeta$,
which sends a tangent vector
$\xi\in T_{0}\R^{n}$
to the tangent vector
$T_{0}T^{r-1}_{n}\zeta\cdot\xi$
to
$T^{r-1}_{n}Y$
at
$J^{r-1}_{0}\zeta$.
If
$\xi =\xi^{i}(\partial /\partial t^{i})_{0}$, then by (2.3),

$$
\gathered
T_{0}T^{r-1}_{n}\zeta\cdot\xi = \sum^{r-1}_{k=0}\, \sum_{i_{1}\le
i_{2}\le\ldots\le i_{k}}\, 
\left({\partial (y^{A}_{i_{1}i_{2}\dots i_{k}}\circ T^{r-1}_{n}\zeta)\over
\partial t^{i}}\right)_{0} \xi^{i} 
\left({\partial \over \partial
y^{A}_{i_{1}i_{2} \dots i_{k}}}\right)_{J^{r-1}_{0}\zeta} \\ \qquad
=\sum^{r-1}_{k=0}\,\sum_{i_{1}\le i_{2}\le\dots\le i_{k}}\,
y^{A}_{i_{1}i_{2}\dots i_{k}i}(J^{r}_{0}\zeta)\xi^{i}
 \left({\partial \over
\partial y^{A}_{i_{1}i_{2}\dots i_{k}}}\right)_{J^{r-1}_{0}\zeta} =
\xi^{i}d_{i}(J^{r}_{0}\zeta),
\endgathered
\tag{(2.4)}
$$
where

$$
d_{i} = \sum^{r-1}_{k=0}\,\sum_{i_{1}\le i_{2} \le \dots\le i_{k}}
y^{A}_{i_{1}i_{2}\dots i_{k}i}
\frac{\partial}{\partial y^{A}_{i_{1}i_{2}\dots i_{k}}} \tag{(2.5)}
$$
is a morphism
$T^{r}_{n}Y \ni J^{r}_{0}\zeta \to d_{i}(J^{r}_{0}\zeta) \in TT^{r-1}_{n}Y$
over
$T^{r-1}_{n}Y$.
Indeed, the tangent vectors
$d_{i}(J^{r}_{0}\zeta)$
are defined independently of the chosen chart: If $(\overline
V,\overline\psi)$,
$\overline\psi =(\overline y^{A})$,
is some other chart at
$\zeta(0)$,
then

$$
\overline d_{i} = \sum^r_{k=0}\,\sum_{j_{1}\le j_{2} \le \dots\le j_{k}}\, 
\overline y^{A}_{j_{1}j_{2}\dots j_{k}i} \frac{\partial}{\partial\overline
y^{A}_{j_{1}j_{2}\dots j_{k}}}, $$
and by (2.4),
$$
\overline d_{i} = d_{i}.
\tag{(2.6)}
$$

We note that formula (2.5) does not define a vector field on $T^{r}_{n}Y$
since it is not invariant when the tangent vectors $\partial /\partial
y^{\sigma}_{i_{1}i_{2}\dots i_{k}}$ are 
subject to coordinate 
transformations on $T^{r}_{n}Y$.

Let
$f: V^{r-1}_{n} \to \R$
be a smooth function. We define the
$i$th
{\it formal derivative}
$d_{i}f: V^{r}_{n} \to \R$
by

$$
d_{i}f = \sum^{r-1}_{k=0} {\sum_{j_{1}\leq j_{2}\leq\ldots\leq j_{k}}
y^{A}_{j_{1} j_{2}\ldots j_{k}i}
{\partial f\over \partial y^{A}_{j_{1}j_{2} \ldots j_{k}}}}. $$

By (2.6), the functions
$d_{i}f$
are independent of the charts, and the definition of the $i$th
formal derivative is naturally extended to functions defined on arbitrary
open subset of
$T^{r-1}_{n}Y$.

It can be easily verified that for every smooth function $f:V^{r-1}_{n} \to
\R$
and every smooth mapping
$\gamma$
of an open set
$U\subset \R^{n}$
into
$Y$,
$
d_{i}f \circ T^{r}_{n}\gamma = D_{i}(f\circ T^{r-1}_{n}\gamma). $
In particular, we have for every coordinate function
$y^{A}_{j_{1}j_{2}\dots j_{k}}$,

$$
d_{i}y^{A}_{j_{1}j_{2}\dots j_{k}} = y^{A}_{j_{1}j_{2}\dots j_{k}i}.
\tag{(2.7)}
$$

Using the formal derivative operators
$d_{i}$,
it is now very easy to find the transformation formulas between two
associated charts on
$T^{r}_{n}Y$.
By (2.7) and (2.6),

$$
\overline y^{A}_{j_{1}j_{2}\dots j_{k}j_{k+1}} = 
\overline d_{j_{k+1}} \overline y^{A}_{j_{1}j_{2}\dots j_{k}} = d_{j_{k+1}}
\overline y^{A}_{j_{1}j_{2}\dots j_{k}} = \dots =d_{j_{k+1}}\dots
d_{j_{2}}d_{j_{1}}\overline y^{A}. $$

This formula may be applied whenever the transformation rules for the
coordinate transformations on
$Y$
are known.

We shall need a formula for higher order partial derivatives of the
composed mapping in a form well 
adapted to its use in various inductive calculations in the higher order
differential geometry and the 
theory of differential invariants.

Let
$n$
and
$k$
be integers. If
$I=\{i_{1},i_{2},\dots,i_{k}\}$
is a set of indices such that
$1\le i_{1},i_{2},\dots,i_{k}\le n$,
we usually denote
$
D_{I}=D_{i_{k}}\dots D_{i_{2}}D_{i_{1}}. $
Let
$U$
and
$V$
be open subsets of
$\R^{n}$,
let
$f:V \to \R$
be a smooth function, and let
$\alpha =(\alpha^{i})$
be a smooth mapping of
$U$
into
$V$.
Then one can prove by induction that

$$
\gathered
D_{i_{s}}\dots D_{i_{2}}D_{i_{1}}(f\circ\alpha)(t) \\ =
\sum^{s}_{k=1}\,\sum_{I=(I_{1},I_{2},\dots,I_{k})}\, 
D_{p_{k}}\dots D_{p_{2}}D_{p_{1}}f(\alpha
(t))D_{I_{k}}\alpha^{p_{k}}(t)\dots D_{I_{2}}\alpha^{p_{2}}(t)D_{I_{1}}
\alpha^{p_{1}}(t), \endgathered
\tag{(2.8)}
$$
where the second sum is extended to all partitions
$(I_{1},I_{2},\dots,I_{k})$
of the set
$\{i_{1},i_{2},\dots,i_{s}\}$.

Let us write the transformation equations from $(V,\psi)$
to
$(\overline V,\overline\psi)$
in the form

$$
\overline y^A=F^A(y^B).
\tag{(2.9)}
$$
We wish to determine explicitly the functions $F_{i_{1}}^{A},
F^{A}_{i_{1}i_{2}},\dots, 
F^{A}_{i_{1}i_{2}\dots i_{r}}$ defining the 
corresponding transformation 

$$
\bar y_{i_{1}i_{2}\ldots i_{k}}^{A} =
F_{i_{1}i_{2}\ldots
i_{k}}^{A}(y^{B},y_{j_{1}}^{B},y_{j_{1}j_{2}}^{B},\ldots,
y_{j_{1}j_{2}\ldots j_{k}}^{B}),
\quad 0\le k\le r \tag{(2.10)}
$$
from
$(V^{r}_{n},\psi^{r}_{n})$
to
$(\overline V^{r}_{n},\overline \psi^{r}_{n})$. 

\proclaim{Lemma 1}
The following formula holds

$$
F_{i_{1}i_{2}\ldots i_{s}}^{A} = \sum\limits_{p=1}^s
{\sum\limits_{(I_{1},I_{2},
\ldots,I_{p})} y_{I_{1}}^{B_{1}}}y_{I_{2}}^{B_{2}}\ldots y_{I_{p}}^{B_{p}}
{{\partial^{p} F^{A}} 
\over {\partial y^{B_{1}}\partial y^{B_{2}}\ldots \partial y^{B_{p}}}},
\tag{(2.11)}
$$
where the second sum denotes summation over all partitions
$(I_{1},I_{2},\dots, I_{p})$
of the set
$\{i_{1},i_{2},\dots,i_{s}\}$.
\endproclaim

\demo{\bf  Proof}
We proceed by induction using (2.7).
\enddemo

We assume that the reader is familiar with the concept of the differential
group. Recall that the 
{\it r\,th differential group of} $\R^{n}$,
denoted by
$L^{r}_{n}$,
is the Lie group of invertible
$r$-jets
with source and target at
$0\in \R^{n}$.
The group multiplication in
$L^{r}_{n}$
is defined by the jet composition

$$
L^{r}_{n}\times L^{r}_{n} \ni (J^{r}_{0}\alpha, J^{r}_{0}\beta) 
\to J^{r}_{0}\alpha\circ J^{r}_{0}\beta = J^{r}_{0}(\alpha\circ \beta) \in
L^{r}_{n}
\tag{(2.12)}
$$
where
$\circ$
denotes both the composition of mappings, and the composition of $r$-
jets. The {\it canonical} (global) {\it coordinates} on $L^{r}_{n}$
are defined by

$$
a^{j}_{i_{1}i_{2}\dots i_{k}}(J^{r}_{0}\alpha) = D_{i_{1}}D_{i_{2}}\dots
D_{i_{k}}\alpha^j(0),\ 1\le k\le r,
\ 1\le i_{1}\le i_{2}\le\dots\le i_{k}\le n, \tag{(2.13)}
$$
where
$\alpha^{j}$
are the components of a representative
$\alpha$
of
$J^{r}_{0}\alpha$.

\proclaim{Lemma 2}
The group multiplication {\rm (2.12)} in
$L^{r}_{n}$
is expressed in the canonical coordinates {\rm (2.13)} by the equations 

$$
c^{k}_{i_{1}i_{2}\dots i_{s}} = \sum^s_{p=1}\,
{\sum_{(I_{1},I_{2},\dots,I_{p})}
\,b^{j_{1}}_{I_{1}}b^{j_{2}}_{I_{2}} \dots b^{j_{p}}_{I_{p}}
a^{k}_{j_{1}j_{2}\dots j_{p}}}, \tag{(2.14)}
$$
where
$a^{k}_{i_{1}i_{2}\dots i_{s}} = a^{k}_{i_{1}i_{2}\dots
i_{s}}(J^{r}_{0}\alpha) $,
$b^{k}_{i_{1}i_{2}\dots i_{s}} = a^{k}_{i_{1}i_{2}\dots
i_{s}}(J^{r}_{0}\beta) $,
$
c^{k}_{i_{1}i_{2}\dots i_{s}} =
a^{k}_{i_{1}i_{2}\dots i_{s}}(J^{r}_{0}(\alpha\circ\beta)) $,
and the second sum is extended to all partitions
$(I_{1},I_{2},\dots,I_{p})$
of the set
$\{i_{1},i_{2},\dots,i_{s}\}$.
\endproclaim

\demo{\bf Proof}
We apply (2.8).
\enddemo

The manifolds of
$(r,n)$-velocities
$T^{r}_{n}Y$
is endowed with a smooth right action of the differential group
$L^{r}_{n}$,
defined by the jet composition

$$
T_{n}^{r}Y\times L_{n}^{r}\,\ni (J_{0}^{r}\zeta,J_{0}^{r}\alpha ) \to
J_{0}^{r}\zeta \circ J_{0}^{r}\alpha = 
J_{0}^{r}(\zeta \circ \alpha ) \in T_{n}^{r}Y.
\tag{(2.15)}
$$
Let us determine the chart expression of this action. To this purpose we
use the canonical coordinates
$a^{i}_{I}$
(2.13) on
$L^{r}_{n}$.

\proclaim{Lemma 3}
The group action {\rm(2.15)} is expressed by the equations 

$$
\bar y^{A} = y^{A},\quad \bar y_{i_{1}i_{2}\ldots i_{s}}^{A} =
\sum\limits_{p=1}^{s} {\sum\limits_{(I_{1},I_{2},
\dots,I_{p})} a_{I_{1}}^{j_{1}}a_{I_{2}}^{j_{2}} \dots a_{I_{p}}^{j_{p}}
y_{j_{1}j_{2}\ldots j_{p}}^{A}},
\tag{(2.16)}
$$
where the second sum is extended to all partitions
$(I_{1},I_{2},\dots,I_{p})$
of the set
$\{i_{1},i_{2},\dots,i_{s}\}$.
\endproclaim

\demo{\bf Proof}
To prove (2.16), we apply (2.1), (2.15) and (2.8). \enddemo

Note the following formula. If
$\gamma $
is a smooth mapping of an open set
$U \subset R^{n}$
into
$Y$,
$U' \subset R^{n}$an open set, and
$\alpha :U'\to U$
a smooth mapping, then for every
$t \in U'$,

$$
T_{n}^{r}(\gamma \circ \alpha )(t) = (T_{n}^{r}\gamma )(\alpha (t)) \circ
J_{0}^{r}(\text{tr}_{-\alpha (t)}\circ 
\alpha \circ\text{tr}_{t}). \tag{(2.17)}
$$

To derive this formula, we use definition (2.2), and the identity $
J_{0}^{r}(\gamma \circ \alpha \circ\text{tr}_{t})\ = J_{0}^{r}(\gamma
\circ\text{tr}_{\alpha (t)})
\circ J_{0}^{r}(\text{tr}_{-\alpha (t)}\circ \alpha \circ\text{tr}_{t}). $
In particular, if
$\alpha$
is a diffeomorphism, then
$
J_{0}^{r}(\text{tr}_{-\alpha (t)} \circ \alpha \circ\text{tr}_{t}) \in
L_{n}^{r}, $
and (2.17) reduces to the group action (2.15). 

\heading 3. Higher order Grassmann bundles \endheading 

An
$(r,n)$-velocity
$J_{0}^{r}\zeta \in T_{n}^{r}Y$
is said to be {\it regular}, if it has a representative which is an
immersion. If
$(V,\psi)$,
$\psi =(y^{A})$,
is a chart, and the target
$\zeta (0)$
of an element
$J_{0}^{r}\zeta \in T_{n}^{r}Y$
belongs to
$V$,
then
$J_{0}^{r}\zeta$
is regular if and only if there exists a subsequence $\nu =(\nu
_{1},\nu_{2},\ldots,\nu_{n})$ of the sequence
$(1,2,\ldots,n,n+1,\ldots,n+m)$
such that
$
\det D_{i}(y^{\nu_{k}}\circ \zeta)(0) \not= 0. $
Regular
$(r,n)$-velocities
form an open,
$L_{n}^{r}$-invariant
subset of
$T_{n}^{r}Y$,
which is called the {\it manifold of regular} $(r,n)$-{\it velocities},
and is denoted by
$\Imm T_{n}^{r}Y$.
Recall that
$\Imm T_{n}^{r}Y$
is endowed with a smooth right action of the differential group
$L_{n}^{r}$,
defined by restricting (2.15), i.e., by

$$
\Imm T_{n}^{r}Y\times L_{n}^{r}\ni (J_{0}^{r}\zeta,J_{0}^{r}\alpha) \to
J_{0}^{r}\zeta \circ J_{0}^{r}\alpha = 
J_{0}^{r}(\zeta \circ \alpha)\in \Imm T_{n}^{r}Y.
\tag{(3.1)}
$$

If
$a_{i_{1}}^{k}, a_{i_{1}i_{2}}^{k},\ldots,a_{i_{1}i_{2}\ldots i_{r}}^{k}$
are the canonical coordinates 
on $L_{n}^{r}$,
this action is expressed by 

$$
\bar y^{A} = y^{A},\quad \bar y_{i_{1}i_{2}\dots i_{s}}^{A} =
\sum\limits_{p=1}^{s} {\sum\limits_{(I_{1},I_{2},
\ldots,I_{p})} a_{I_{1}}^{j_{1}}a_{I_{2}}^{j_{2}}\dots a_{I_{p}}^{j_{p}}
y_{j_{1}j_{2}\dots j_{p}}^{A}}.
\tag{(3.2)}
$$

In the proof of the following result we construct, among others, a complete
system of
$L_{n}^{r}$-invariants
of the action (3.1) on
$\Imm T_{n}^{r}Y$.
We use the {\it associated charts} on
$\Imm T_{n}^{r}Y$, 
defined as associated charts on
$T_{n}^{r}Y$.

\proclaim{Theorem 1}
The group action {\rm(3.1)} defines on
$\Imm T_{n}^{r}Y$
the structure of a right principal
$L_{n}^{r}$-bundle.
\endproclaim

\demo{\bf Proof}
We have to show that (a) the equivalence $\Cal R$
``there exists
$J_{0}^{r}\alpha \in L_{n}^{r}$
such that
$J_{0}^{r}\zeta = J_{0}^{r}\chi \circ J_{0}^{r}\alpha$ is a closed
submanifold of the product manifold 
$\Imm T_{n}^{r}Y\times \Imm T_{n}^{r}Y$, and (b) the group action (3.1) is
free.

\paritem[(a)] First we construct an atlas on $\Imm T_{n}^{r}Y\times \Imm
T_{n}^{r}Y$, adapted to the group action 
(3.1).

Let
$(V,\psi)$, $\psi = (y^{A})$,
be a chart on
$Y$,
$(V_{n}^{r},\psi_{n}^{r})$,
$\psi^{r}_{n} = (y^{A},y^{A}_{j_{1}},\ldots,y^{A}_{j_{1}j_{2}\ldots
j_{r}})$, the associated chart on
$\Imm T_{n}^{r}Y$.
We set for every subsequence
$\nu =(\nu_{1},\nu_{2},\ldots,\nu_{n})$
of the sequence
$(1,2,\ldots,n,n+1,\ldots,n+m)$

$$
W^{\nu} = \left\{J_{0}^{r}\zeta \in V_{n}^{r} |
\det(y_{j}^{\nu_{k}}(J_{0}^{r}\zeta) \not= 0 \right\}. \tag{(3.3)}
$$
$W^{\nu}$
is an open,
$L_{n}^{r}$-invariant
subset of
$V_{n}^{r}$,
and

$$
\bigcup\limits_{\nu} {W^{\nu}} = V_{n}^{r}. $$
Restricting the mapping
$\psi_{n}^{r}$
to
$W^{\nu}$
we obtain a chart
$(W^{\nu},\psi_{n}^{r})$.

The equivalence
$\Cal R$
is obviously covered by the open sets of the form $W^{\nu} \times W^{\nu}$.
We shall find its equations in terms of the charts $(W^{\nu} \times
W^{\nu},\psi_{n}^{r}\times \psi_{n}^{r})$. 
Let us consider the set
$\Cal R \cap (W^{\nu} \times W^{\nu})$. Assume for simplicity that
$\nu =(1,2,\ldots,n)$.
A point
$(J_{0}^{r}\zeta,J_{0}^{r}\chi )\in W^{\nu} \times W^{\nu}$ belongs to
$\Cal R \cap (W^{\nu} \times W^{\nu})$
if and only if there exists
$J_{0}^{r}\alpha \in L_{n}^{r}$
such that
$J_{0}^{r}\zeta = J_{0}^{r}\chi \circ J_{0}^{r}\alpha$ or, which is the
same, if and only if the system of 
equations (3.2) has a solution
$a_{i_{1}}^{k},a_{i_{1}i_{2}}^{k},\ldots,a_{i_{1}i_{2}\ldots i_{r}}^{k}$.
Clearly, in this system
$
\bar y^{A},\bar y_{i_{1}}^{A},\bar y_{i_{1}i_{2}}^{A},\ldots, \bar
y_{i_{1}i_{2}\ldots i_{r}}^{A}
$
(resp.
$y^{A},y_{p_{1}}^{A},y_{p_{1}p_{2}}^{A},\ldots,y_{p_{1}p_{2}\ldots
p_{r}}^{A} $)
are coordinates of
$J_{0}^{r}\zeta$
(resp.
$J_{0}^{r}\chi$).
But on
$W^{\nu}$,
$\det(y_{i}^{k}) \not= 0,$
where
$1\le i,k\le n.$
Consequently, there exist functions
$z_{j}^{i}:W^{\nu} \to \R$
such that
$
z_{j}^{i} y_{i}^{k} = \delta_{j}^{k}.
$
Conditions (3.2) now imply, for
$A = k = 1,2,\ldots,n$,

$$
\gathered
\bar y_{i_{1}i_{2}\ldots i_{s}}^{k} = \sum\limits_{p=1}^{s}
{\sum\limits_{(I_{1}I_{2},
\ldots,I_{p})} a_{I_{1}}^{j_{1}}a_{I_{2}}^{j_{2}} \ldots a_{I_{p}}^{j_{p}}
y_{j_{1}j_{2}\ldots j_{p}}^{k}} \\ =
\sum\limits_{p=2}^{s} {\sum\limits_{(I_{1}I_{2},
\dots,I_{p})} a_{I_{1}}^{j_{1}}a_{I_{2}}^{j_{2}}\dots a_{I_{p}}^{j_{p}}
y_{j_{1}j_{2}\ldots j_{p}}^{k}}
+ a_{i_{1}i_{2} \dots i_{s}}^{j_{1}} y_{j_{1}}^{k}, \endgathered
$$
which allows to determine the canonical coordinates of the group element
$J_{0}^{r}\alpha$
by the recurrent formula

$$
a_{i_{1}i_{2}\ldots i_{s}}^{q} = z_{k}^{q} \Big( \bar y_{i_{1}i_{2}\ldots
i_{s}}^{k}
-\sum\limits_{p=2}^{s} {\sum\limits_{(I_{1}I_{2},\ldots,I_{p})}
a_{I_{1}}^{j_{1}}a_{I_{2}}^{j_{2}}
\ldots a_{I_{p}}^{j_{p}} y_{j_{1}j_{2}\ldots j_{p}}^{k}}\Big).
\tag{(3.4)}
$$

Taking
$A = \sigma = n+1, n+2, \ldots, n+m$
in (3.2) and substituting from (3.4) we get 

$$
\bar y_{i_{1}i_{2}\ldots i_{s}}^{\sigma} = \sum\limits_{p=1}^{s}
\sum\limits_{(I_{1},I_{2},
\ldots,I_{p})} a_{I_{1}}^{j_{1}}a_{I_{2}}^{j_{2}} \ldots a_{I_{p}}^{j_{p}}
y_{j_{1}j_{2}\ldots j_{p}}^{\sigma}, 
\tag{(3.5)}
$$
where the group parameters
$a_{I}^{j}$
are all certain rational functions of
$
y_{j_{1}j_{2}\ldots j_{s}}^{\lambda},
\bar y_{j_{1}j_{2}\ldots j_{s}}^{\lambda}$. These are the desired equations
of the equivalence $\Cal R$
on
$W^{\nu} \times W^{\nu}$.

Now define a new chart on
$\Imm T_{n}^{r}Y\times \Imm T_{n}^{r}Y$, $(W^{\nu} \times
W^{\nu},\Phi^{\nu})$,
where

$$
\Phi^{\nu} = (y^{A},y_{j_{1}}^{A},y_{j_{1}j_{2}}^{A},\ldots,
y_{j_{1}j_{2}\ldots j_{r}}^{A},
\Phi^\sigma,\Phi_{j_{1}}^{\sigma},\Phi_{j_{1}j_{2}}^{\sigma},\ldots,
\Phi_{j_{1}j_{2}\ldots j_{r}}^{\sigma},
\bar y^{k},\bar y_{j_{1}}^{k},\bar y_{j_{1}j_{2}}^{k},\ldots, \bar
y_{j_{1}j_{2}\ldots j_{r}}^{k}),
$$
is the collection of coordinate functions, by 

$$
\Phi^{\sigma} = \bar y^{\sigma} - y^{\sigma},\quad \Phi_{i_{1}i_{2}\ldots
i_{s}}^{\sigma} = 
\bar y_{i_{1}i_{2}\ldots i_{s}}^{\sigma} - \sum\limits_{p=1}^{s}
{\sum\limits_{(I_{1},I_{2},
\ldots,I_{p})} a_{I_{1}}^{j_{1}}a_{I_{2}}^{j_{2}} \ldots a_{I_{p}}^{j_{p}}
y_{j_{1}j_{2}\ldots j_{p}}^{\sigma}}. $$
In terms of this new chart, the equivalence $\Cal R$
has equations
$
\Phi^\sigma = 0,\Phi_{i_{1}i_{2}\ldots i_{s}}^{\sigma} = 0, $
and is therefore a closed submanifold of $\Imm T_{n}^{r}Y\times \Imm
T_{n}^{r}Y$. 

\paritem[(b)] Assume that for some
$J_{0}^{r}\zeta \in \Imm T_{n}^{r}Y$
and
$
J_{0}^{r}\alpha \in L_{n}^{r}, J_{0}^{r}\zeta \circ J_{0}^{r}\alpha =
J_{0}^{r}\zeta.
$
Then equations (3.2) reduce to

$$
y_{i_{1}i_{2}\ldots i_{s}}^{A} = \sum\limits_{p=1}^{s}
{\sum\limits_{(I_{1},I_{2},
\ldots,I_{p})} a_{I_{1}}^{j_{1}}a_{I_{2}}^{j_{2}}\ldots a_{I_{p}}^{j_{p}}
y_{j_{1}j_{2}\ldots j_{p}}^{A}},
$$
which gives us, using (3.4),
$
a_{i}^{p} = \delta_{i}^{p}, a_{i_{1}i_{2}}^{p} = 0, \ldots,
a_{i_{1}i_{2}\ldots i_{r}}^{p} = 0$,
i.e.,
$J_{0}^{r}\alpha = J_{0}^{r}\id_{\R^{n}}$. This completes the proof.
\enddemo

We have the following corollaries.

\proclaim{Corollary 1}
The orbit space
$P_{n}^{r}Y = \Imm T_{n}^{r}Y/L_{n}^{r}$ has a unique smooth structure such
that the canonical quotient 
projection $\rho_{n}^{r}: \Imm T_{n}^{r}Y \to P_{n}^{r}Y$ is a surjective
submersion. 
The dimension of $P_{n}^{r}Y$
is

$$
\dim P_{n}^{r}Y = m{\pmatrix {n+r}\\ n\endpmatrix}+n. $$
\endproclaim

The following corollary solves the problem of finding all $L_{n}^{r}$-{\it
invariant functions} on 
$\Imm T_{n}^{r}Y$.
It says that the projection
$\rho^{r}_{n}: \Imm T_{n}^{r} \to P_{n}^{r}Y$ is the {\it basis} of
$L_{n}^{r}$-{\it invariant functions}.

\proclaim{Corollary 2}
Every
$L_{n}^{r}$-invariant
function
$f:W \to \R$,
where
$W \subset \Imm T_{n}^{r}Y$
is an
$L_{n}^{r}$-invariant
open set, can be factored through the projection mapping $\rho_{n}^{r}:
\Imm T_{n}^{r}Y \to P_{n}^{r}Y$. 
\endproclaim

Now we are going to construct charts on
$\Imm T_{n}^{r}Y$
adapted to the right action (3.1) of the differential group $L_{n}^{r}$.
We may consider, for example, the charts (3.3) with $\nu =(1,2,\ldots,n)$.

\proclaim{Theorem 2}
Let
$(V,\psi), \psi = (y^{A})$,
be a chart on
$
Y, (V_{n}^{r},\psi_{n}^{r}), \psi_{n}^{r} = (y^{A}_{i_{1}i_{2}\dots
i_{s}}), s \le r,
$
the associated chart on
$\Imm T_{n}^{r}Y$,
and
$
W = \left\{ J_{0}^{r}\zeta \in V_{n}^{r}| \det (y_{j}^{i}(J_{0}^{r}\zeta))
\not= 0 \right\},1\le i,j\le n.
$
There exist unique functions
$
w^{\sigma}, w_{j_{1}}^{\sigma}, w_{j_{1}j_{2}}^{\sigma}, \ldots,
w_{j_{1}j_{2}\ldots j_{r}}^{\sigma} $
defined on
$W$
such that

$$
y^{\sigma} = w^{\sigma},\quad y_{p_{1}p_{2}\ldots p_{k}}^{\sigma} =
\sum\limits_{q=1}^{k} 
\sum\limits_{(I_{1},I_{2},\ldots,I_{q})}
y_{I_{1}}^{j_{1}}y_{I_{2}}^{j_{2}}\ldots y_{I_{q}}^{j_{q}} 
w_{j_{1}j_{2}\ldots j_{q}}^{\sigma}.
\tag{(3.6)}
$$
The pair
$(W,\Phi)$,
where
$
\Phi = (w^{\sigma}, w_{p_{1}}^{\sigma},w_{p_{1}p_{2}}^{\sigma},\ldots,
w_{p_{1}p_{2}\ldots p_{r}}^{\sigma},
y^{i},y_{j_{1}}^{i},y_{j_{1}j_{2}}^{i}, \ldots,y_{j_{1}j_{2}\ldots
j_{r}}^{i})
$,
is a chart on  \linebreak
$\Imm T_{n}^{r}Y$.
The functions
$
w^{\sigma}, w_{j_{1}}^{\sigma},w_{j_{1}j_{2}}^{\sigma},\ldots,
w_{j_{1}j_{2}\ldots j_{r}}^{\sigma}
$
satisfy the recurrent formula

$$
w_{j_{1}j_{2}\ldots j_{k}j_{k+1}}^{\sigma} = z_{j_{k+1}}^{s} d_{s}
w_{j_{1}j_{2}\ldots j_{k}}^{\sigma},
\tag{(3.7)}
$$
and are
$L_{n}^{r}$-invariant.
\endproclaim

\demo{\bf Proof} We proceed by induction.

\paritem[1.] We prove that the assertion is true for $r=1$.
Consider the pair
$(W,\Phi), \Phi = (w^{\sigma},w_{p_{1}}^{\sigma},y^{i},y_{j_{1}}^{i})$,
where
$w^{\sigma} =y^{\sigma}, w_{j}^{\sigma} = z_{j}^{k} y_{k}^{\sigma}$.
Obviously
$y_{p}^{\sigma} = y_{p}^{j}w_{j}^{\sigma}$, which implies that
$(W,\Phi)$
is a new chart. Moreover,
$w_{j}^{\sigma} = z_{j}^{k}d_{k}y^{\sigma} = z_{j}^{k}d_{k}w^{\sigma}$. It
remains to show that the functions
$w_{j}^{\sigma}$
are
$L_{n}^{1}$-invariant.
Since the group action (3.2) is expressed by $
\bar y^{i} = y^{i}, \bar y^{\sigma} = y^{\sigma}, \bar y_{p}^{i} =
a_{p}^{j}y_{j}^{i},
\bar y_{p}^{\sigma} =a_{p}^{j}y_{j}^{\sigma}$, the inverse of the matrix
$\bar y_{p}^{i} = a_{p}^{j}y_{j}^{i}$
is
$\bar z_{q}^{p} = z_{q}^{s}b_{s}^{p}$,
where
$b_{s}^{p}$
stands for the inverse of
$a_{s}^{p}$.
Hence
$
\bar w_{j}^{\sigma} = \bar z_{j}^{k}\bar y_{k}^{\sigma} =
z_{j}^{s}b_{s}^{k}a_{k}^{p}y_{p}^{\sigma} = 
z_{j}^{p}y_{p}^{\sigma} = w_{j}^{\sigma}
$
proving the invariance.

\paritem[2.] Assume that formulas (3.6), (3.7) hold for $k=r-1$.
Write (3.6) in the form

$$
y_{p_{1}p_{2}\ldots p_{k}}^{\sigma} = \sum\limits_{q=1}^{k}
\sum\limits_{(I_{1},I_{2},\ldots,
I_{q})} y_{I_{1}}^{j_{1}}y_{I_{2}}^{j_{2}} \ldots y_{I_{q}}^{j_{q}}
w_{j_{1}j_{2}\ldots j_{q}}^{\sigma}. $$

Then

$$
\gathered
y_{p_{1}p_{2}\ldots p_{k}p_{k+1}}^{\sigma} = d_{p_{k+1}}
y_{p_{1}p_{2}\ldots p_{k}}^{\sigma} \\
= \sum\limits_{q=1}^{k} {\sum\limits_{(I_{1},I_{2},\ldots,I_{q})}
(d_{p_{k+1}} 
(y_{I_{1}}^{j_{1}}y_{I_{2}}^{j_{2}}\ldots y_{I_{q}}^{j_{q}})
w_{j_{1}j_{2}\ldots j_{q}}^{\sigma} +
y_{I_{1}}^{j_{1}}y_{I_{2}}^{j_{2}}\ldots y_{I_{q}}^{j_{q}}
y_{p_{k+1}}^{j_{q+1}} z_{j_{q+1}}^{s} 
d_{s} w_{j_{1}j_{2}\ldots j_{q}}^{\sigma})}.
\endgathered
$$

In this formula we sum through all partitions $(I_{1},I_{2},\ldots,I_{q})$
of the set
$\{ p_{1},p_{2},\ldots,p_{k}\}.$
On the other hand, when passing to all partitions
$(J_{1},J_{2},\ldots,J_{q})$
of the set
$\{ p_{1},p_{2},\ldots,p_{k},p_{k+1}\}$
we get

$$
\gathered
y^{\sigma}_{p_{1}p_{2}\ldots p_{k}p_{k+1}} \\ = \sum\limits_{q=1}^{k}
{\sum\limits_{(I_{1},I_{2},\ldots,I_{q})} 
(d_{p_{k+1}} (y_{I_{1}}^{j_{1}}y_{I_{2}}^{j_{2}}\ldots y_{I_{q}}^{j_{q}})
w_{j_{1}j_{2}\ldots j_{q}}^{\sigma} +
y_{I_{1}}^{j_{1}}y_{I_{2}}^{j_{2}}
\ldots y_{I_{q}}^{j_{q}} y_{p_{k+1}}^{j_{q+1}}z_{j_{q+1}}^{s} d_{s}
w_{j_{1}j_{2}\ldots j_{q}}^{\sigma})} \\ = 
\sum\limits_{q=1}^{k} {\sum\limits_{(J_{1},J_{2},\ldots,J_{q})}
y_{J_{1}}^{j_{1}}y_{J_{2}}^{j_{2}}
\ldots y_{J_{q}}^{j_{q}} w_{j_{1}j_{2}\ldots j_{q}}^{\sigma}} +
y_{p_{1}}^{j_{1}}y_{p_{2}}^{j_{2}} 
\ldots y_{p_{q}}^{j_{q}} y^{j_{q+1}}_{p_{k+1}} z_{j_{q+1}}^{s} d_{s}
w_{j_{1}j_{2}\ldots j_{q}}^{\sigma},
\endgathered
\tag{(3.8)}
$$
and we see that (3.8) has the same form as (3.6), where $
w_{j_{1}j_{2}\ldots j_{k}j_{k+1}}^{\sigma} = z_{j_{k+1}}^{s}d_{s}
w_{j_{1}j_{2}\ldots j_{k}}^{\sigma}.
$
Uniqueness is immediate since
$w_{j_{1}j_{2}\ldots j_{k}j_{k+1}}^{\sigma}$ may be expressed explicitly
from (3.8).

It remains to prove the invariance condition $
\bar w_{j_{1}j_{2}\ldots j_{s}}^{\sigma} = w_{j_{1}j_{2}\ldots
j_{s}}^{\sigma}
$.

Since the points
$J^{r}_{0}\zeta\circ J^{r}_{0}\alpha$
and
$J^{r}_{0}\zeta$
belong to the same orbit, their coordinates satisfy the recurrence formula
(3.5):

$$
\bar y^{\sigma}_{i_{1}i_{2}\ldots i_{s}} = \sum\limits_{p=1}^{s}
{\sum\limits_{(I_{1},I_{2},\ldots,I_{p})} 
a_{I_{1}}^{j_{1}}a_{I_{2}}^{j_{2}}\ldots a_{I_{p}}^{j_{p}}
y_{j_{1}j_{2}\ldots j_{p}}^{\sigma}}, \quad s = 
1,2,\ldots,r, $$
in which

$$
a^{q}_{k_{1}k_{2}\ldots k_{t}} = z^{q}_{k} \Big( \bar
y^{k}_{i_{1}i_{2}\ldots i_{t}} - \sum\limits_{p=2}^{t} 
{\sum\limits_{(K_{1},K_{2},\ldots,K_{p})}
a_{I_{1}}^{j_{1}}a_{I_{2}}^{j_{2}}\ldots a_{I_{p}}^{j_{p}} y_{j_{1}j_{2}
\ldots j_{p}}^{k}} \Big)
$$
for all
$ t \le s$
(see (3.4)). Here
$(I_{1},I_{2},\ldots I_{p})$
is a partition of the set
$\{ i_{1},i_{2},\ldots i_{s}\}$
and
$(K_{1},K_{2},\ldots K_{p})$
$\{ k_{1},k_{2},\ldots k_{t}\}.$
Using (3.6) we can write

$$
\gathered
\bar y_{i_{1}i_{2}\ldots i_{s}}^{\sigma} = \sum\limits_{p=1}^{s}
{\sum\limits_{(I_{1},I_{2},\ldots,I_{p})} 
\bar y_{I_{1}}^{j_{1}}\bar y_{I_{2}}^{j_{2}}\ldots \bar y_{I_{p}}^{j_{p}}
\bar w_{j_{1}j_{2}
\ldots j_{p}}^{\sigma}}, \\ y_{j_{1}j_{2}\ldots j_{s}}^{\sigma} =
\sum\limits_{l=1}^{p} {\sum\limits_{(J_{1},J_{2},\ldots,J_{l})}
y_{J_{1}}^{t_{1}}y_{J_{2}}^{t_{2}}
\ldots y_{J_{l}}^{t_{l}} w_{t_{1}t_{2}\ldots t_{l}}^{\sigma}} ,
\endgathered
$$
where
$(I_{1},I_{2},\ldots I_{p})$
is a partition of the set
$\{ i_{1},i_{2},\ldots i_{s}\}$
and
$(J_{1},J_{2},\ldots J_{l})$
is a partition of the set
$\{ j_{1},j_{2},\ldots j_{p}\}.$
This gives us the equation

$$
\gathered
\sum\limits_{p=1}^{s} {\sum\limits_{(I_{1},I_{2},\ldots,I_{p})} \bar
y_{I_{1}}^{j_{1}}\bar y_{I_{2}}^{j_{2}}
\ldots \bar y_{I_{p}}^{j_{p}} \bar w_{j_{1}j_{2}\ldots j_{p}}^{\sigma}} \\
= \sum\limits_{p=1}^{s} 
{\sum\limits_{(I_{1},I_{2},\ldots,I_{p})}
a_{I_{1}}^{j_{1}}a_{I_{2}}^{j_{2}}\ldots a_{I_{p}}^{j_{p}}} 
\Big( \sum\limits_{l=1}^{p} {\sum\limits_{(J_{1},J_{2},\ldots,J_{l})}
y_{J_{1}}^{t_{1}}y_{J_{2}}^{t_{2}}
\ldots y_{J_{l}}^{t_{l}} w_{t_{1}t_{2}\ldots t_{l}}^{\sigma}} \Big).
\endgathered
\tag{(3.9)}
$$

Now we wish to determine the terms
$w_{t_{1}t_{2}\ldots t_{p}}^{\sigma}$
on the right side with fixed
$p$.
Changing the notation of the indices, we get the expression 

$$
\sum\limits_{q=1}^{s} {\sum\limits_{(I_{1},I_{2},\ldots,I_{p})}
a_{I_{1}}^{j_{1}}a_{I_{2}}^{j_{2}}
\ldots a_{I_{q}}^{j_{q}}} \Big( \sum\limits_{p=1}^{q}
{\sum\limits_{(J_{1},J_{2},
\ldots,J_{p})} y_{J_{1}}^{t_{1}}y_{J_{2}}^{t_{2}}\ldots y_{J_{p}}^{t_{p}}
w_{t_{1}t_{2}\ldots t_{p}}^{\sigma}} 
\Big) $$
from which we see that
$w_{t_{1}t_{2}\ldots t_{p}}^{\sigma}$
are contained in every summand with
$ q \ge p.$
Thus, the required terms are given by

$$
\Big( \sum\limits_{q=p}^{s} \sum\limits_{(I_{1},I_{2},\ldots,I_{q})}
{\sum\limits_{(J_{1},J_{2},\ldots,J_{p})} 
a_{I_{1}}^{j_{1}}a_{I_{2}}^{j_{2}}\ldots a_{I_{q}}^{j_{q}}
y_{J_{1}}^{t_{1}}y_{J_{2}}^{t_{2}}
\ldots y_{J_{p}}^{t_{p}}}\Big) w_{t_{1}t_{2}\ldots t_{p}}^{\sigma}.
$$

In this formula
$(I_{1},I_{2},\ldots I_{q})$
is a partition of the set
$\{ i_{1},i_{2},\ldots i_{s}\},$
and
$(J_{1},J_{2},\ldots J_{p})$
is a partition of the set
$\{ j_{1},j_{2},\ldots j_{q}\}.$

Now we adopt the following notation. If
$I = (i_{1},i_{2},\ldots i_{s})$
is a multi-index, then 
$(I_{1},I_{2},\ldots I_{p}) \sim I$
means that
$(I_{1},I_{2},\ldots I_{p})$
is a partition of the set
$\{ i_{1},i_{2},\ldots i_{s}\}.$

As before, let
$I = (i_{1},i_{2},\ldots i_{s})$,
and let
$p$
be fixed. We wish to show that

$$
\gathered
\Big( {\sum\limits_{(I_{1},I_{2},\ldots,I_{p})} \bar y_{I_{1}}^{j_{1}}\bar
y_{I_{2}}^{j_{2}}\ldots 
\bar y_{I_{p}}^{j_{p}}} \Big) w_{j_{1}j_{2}\ldots j_{p}}^{\sigma} \\  =
\Big( \sum\limits_{q=p}^{s} 
{\sum\limits_{(I_{1},I_{2},\ldots,I_{q})}} {\sum\limits_{(J_{1},J_{2},
\ldots,J_{p})} a_{I_{1}}^{j_{1}}a_{I_{2}}^{j_{2}}\ldots a_{I_{q}}^{j_{q}}
y_{J_{1}}^{t_{1}}y_{J_{2}}^{t_{2}}
\ldots y_{J_{p}}^{t_{p}}}\Big) w_{t_{1}t_{2}\ldots t_{p}}^{\sigma}.
\endgathered
\tag{(3.10)}
$$

Write the transformation formula (2.16) in the form 

$$
\bar y_{I}^{A} =
\sum\limits_{p=1}^{|I|} {\sum\limits_{(I_{1},I_{2},\ldots,I_{p})}
a_{I_{1}}^{j_{1}}a_{I_{2}}^{j_{2}} 
\ldots a_{I_{p}}^{j_{p}} y_{j_{1}j_{2}\ldots j_{p}}^{A}}, \quad
(I_{1},I_{2},\ldots I_{p}) \sim I.
$$
Using the same notation, we have

$$
\bar y_{I_{k}}^{t_{k}}
=
\sum\limits_{q_{k}=1}^{|I_{k}|}{\sum\limits_{(I_{k,1},I_{k,2},\ldots,I_{k,q_
{k}})} 
a_{I_{k,1}}^{j_{k,1}}a_{I_{k,2}}^{j_{k,2}} \ldots
a_{I_{k,q_{k}}}^{j_{k,q_{k}}}
y_{j_{k,1}j_{k,2}\ldots j_{k,q_{k}}}^{t_{k}}}, \quad
(I_{k,1},I_{k,2},\ldots I_{k,q_{k}}) \sim I_{k}, $$
where
$(I_{1},I_{2},\ldots I_{p}) \sim I.$
Thus,

$$
\gathered
\Big( {\sum\limits_{(I_{1},I_{2},\ldots,I_{p})} \bar y_{I_{1}}^{j_{1}}\bar
y_{I_{2}}^{j_{2}}\ldots 
\bar y_{I_{p}}^{j_{p}}} \Big) w_{j_{1}j_{2}\ldots j_{p}}^{\sigma} = \Big(
\sum\limits_{q_{1}=1}^{|I_{1}|}
{\sum\limits_{(I_{1,1},I_{1,2},\ldots,I_{1,q_{1}})}
a_{I_{1,1}}^{j_{1,1}}a_{I_{1,2}}^{j_{1,2}}
\ldots a_{I_{1,q_{1}}}^{j_{1,q_{1}}} y_{J_{1}}^{t_{1}}}\Big)\\
\Big( \sum\limits_{q_{2}=1}^{|I_{2}|}
{\sum\limits_{(I_{2,1},I_{2,2},\ldots,I_{2,q_{2}})}
a_{I_{2,1}}^{j_{2,1}}a_{I_{2,2}}^{j_{2,2}}
\ldots a_{I_{2,q_{2}}}^{j_{2,q_{2}}} y_{J_{2}}^{t_{2}}}\Big)\\
\ldots \Big( \sum\limits_{q_{p}=1}^{|I_{p}|}
{\sum\limits_{(I_{p,1},I_{p,2},\ldots,I_{p,q_{p}})} 
a_{I_{p,1}}^{j_{p,1}}a_{I_{p,2}}^{j_{p,2}}\ldots
a_{I_{p,q_{p}}}^{j_{p,q_{p}}} y_{J_{p}}^{t_{p}}}\Big)
w_{t_{1}t_{2}\ldots t_{p}}^{\sigma},
\endgathered
$$
where
$
J_{1} = (j_{1,1}j_{1,2}\ldots j_{1,q_{1}}), J_{2} = (j_{2,1}j_{2,2}\ldots
j_{2,q_{2}}),\ldots $
and
$J_{p} = (j_{p,1}j_{p,2}\ldots j_{p,q_{p}}).$ 

This expression can be written in a different way. Notice that since
$(I_{i,1},I_{i,2},\ldots I_{i,q_{i}}) 
\sim I_{i}$, then

$$
\gathered
(I_{1,1},I_{1,2},\ldots I_{1,q_{1}},I_{2,1},I_{2,2},\ldots I_{2,q_{2}},
\ldots I_{p,1},I_{p,2},\ldots I_{p,q_{p}}) \sim I \\ |I_{1}| + |I_{2}| +
\ldots + |I_{p}| = |I| = s. \endgathered
$$
and if we define
$q = q_{1} + q_{2} +\ldots + q_{p}$,
we get
$p \le q \le |I_{1}| + |I_{2}| + \ldots + |I_{p}| = |I| = s.$ Now, having
in mind the corresponding summation 
ranges, 

$$
\gathered
\Big( {\sum\limits_{(I_{1},I_{2},\ldots,I_{p})} \bar y_{I_{1}}^{j_{1}}\bar
y_{I_{2}}^{j_{2}}\ldots 
\bar y_{I_{p}}^{j_{p}}} \Big) w_{j_{1}j_{2}\ldots j_{p}}^{\sigma} \\ = \sum
a_{I_{1,1}}^{j_{1,1}}a_{I_{1,2}}^{j_{1,2}}\ldots
a_{I_{1,q_{1}}}^{j_{1,q_{1}}} 
a_{I_{2,1}}^{j_{2,1}}a_{I_{2,2}}^{j_{2,2}}\ldots
a_{I_{2,q_{2}}}^{j_{2,q_{2}}} \ldots
a_{I_{p,1}}^{j_{p,1}}a_{I_{p,2}}^{j_{p,2}}\ldots
a_{I_{p,q_{p}}}^{j_{p,q_{p}}} 
y_{J_{1}}^{t_{1}}y_{J_{2}}^{t_{2}}\ldots y_{J_{p}}^{t_{p}}
w_{t_{1}t_{2}\ldots t_{p}}^{\sigma}.
\endgathered
$$

If we denote
$$
(s_{1},s_{2},\ldots s_{q}) =
(j_{1,1},j_{1,2},\ldots j_{1,q_{1}},j_{2,1},j_{2,2},\ldots j_{2,q_{2}},
\ldots j_{p,1},j_{p,2},
\ldots j_{p,q_{p}}) $$
and
$$
(P_{1},P_{2},\ldots P_{q}) =
(I_{1,1},I_{1,2},\ldots I_{1,q_{1}},I_{2,1},I_{2,2},\ldots I_{2,q_{2}},
\ldots I_{p,1},I_{p,2},
\ldots I_{p,q_{p}}) $$
it is immediate that
$(P_{1},P_{2},\ldots P_{q}) \sim I$,
and

$$
\gathered
\Big( {\sum\limits_{(I_{1},I_{2},\ldots,I_{p})} \bar y_{I_{1}}^{j_{1}}\bar
y_{I_{2}}^{j_{2}}\ldots 
\bar y_{I_{p}}^{j_{p}}} \Big) w_{j_{1}j_{2}\ldots j_{p}}^{\sigma} \\ =
\Big( \sum\limits_{q=p}^{s}
{\sum_{(P_{1},P_{2},\ldots P_{q})}} {\sum_{(J_{1},J_{2},\ldots J_{p})}
a_{P_{1}}^{s_{1}}a_{P_{2}}^{s_{2}}
\ldots a_{P_{q}}^{s_{q}} y_{J_{1}}^{t_{1}}y_{J_{2}}^{t_{2}}\ldots
y_{J_{p}}^{t_{p}}}\Big) w_{t_{1}t_{2}
\ldots t_{p}}^{\sigma}.
\endgathered
$$
This proves (3.10).

Returning to (3.9), and substituting from (3.10) we get a basic formula 

$$
\sum\limits_{p=1}^{s}
{\sum_{(I_{1},I_{2},\ldots I_{p})}
\bar y_{I_{1}}^{j_{1}}\bar y_{I_{2}}^{j_{2}}\ldots \bar y_{I_{p}}^{j_{p}}
(\bar w^{\sigma}_{j_{1}j_{2}
\ldots j_{p}} - w^{\sigma}_{j_{1}j_{2}\ldots j_{p}}) = 0}. \tag{(3.11)}
$$

Now it is easy to show that
$
\bar w^{\sigma}_{j_{1}j_{2}\ldots j_{s}} - w^{\sigma}_{j_{1}j_{2}\ldots
j_{s}} = 0
$
provided
$
\bar w^{\sigma}_{j_{1}j_{2}\ldots j_{k}} - w^{\sigma}_{j_{1}j_{2}\ldots
j_{k}} = 0
$
for all
$k \le s-1.$

If
$s = 1$,
we get
$
\bar y_{i_{1}}^{j_{1}} (\bar w^{\sigma}_{j_{1}} - w^{\sigma}_{j_{1}}) = 0,
$
and since the matrix
$\bar y_{i}^{j}$
is regular,
$
\bar w^{\sigma}_{j} = w^{\sigma}_{j}.
$

If
$s = 2$,
we have
$
\bar y_{i_{1}i_{2}}^{j_{1}} (\bar w^{\sigma}_{j_{1}} - w^{\sigma}_{j_{1}})
+
\bar y_{i_{1}}^{j_{1}} \bar y_{i_{2}}^{j_{2}} (\bar w^{\sigma}_{j_{1}j_{2}}
- w^{\sigma}_{j_{1}j_{2}}) =
\bar y_{i_{1}}^{j_{1}} \bar y_{i_{2}}^{j_{2}} (\bar w^{\sigma}_{j_{1}j_{2}}
- w^{\sigma}_{j_{1}j_{2}}) = 0,
$
which implies, again using regularity of the matrix $\bar y_{i}^{j}$,
that
$\bar w^{\sigma}_{j_{1}j_{2}} = w^{\sigma}_{j_{1}j_{2}}.$ 

Now assume that
$
\bar w^{\sigma}_{j_{1}j_{2}\ldots j_{k}} - w^{\sigma}_{j_{1}j_{2}\ldots
j_{k}} = 0
$
for all
$k \le s-1.$
Then (3.11) reduces to

$$
\bar y_{i_{1}}^{j_{1}}\bar y_{i_{2}}^{j_{2}}\ldots \bar y_{i_{s}}^{j_{s}}
(\bar w^{\sigma}_{j_{1}j_{2}
\ldots j_{s}} - w^{\sigma}_{j_{1}j_{2}\ldots j_{s}}) = 0, $$
which gives us
$
\bar w^{\sigma}_{j_{1}j_{2}\ldots j_{s}} - w^{\sigma}_{j_{1}j_{2}\ldots
j_{s}} = 0
$
as required.

This completes the proof.
\enddemo

Denote

$$
\Delta_{i} = z^{s}_{i}d_{s}.
\tag{(3.12)}
$$
Properties of the group action of
$L_{n}^{r}$
on
$\Imm T_{n}^{r}Y$
can now be summarized as follows.

\proclaim{Corollary 1}
The group action {\rm(3.1)} is expressed on
$W$
by the equations

$$
\gathered
\bar y^{k}=y^{k},\quad \bar y_{i_{1}i_{2}\ldots i_{s}}^{k}=
\sum\limits_{p=1}^{s} {\sum\limits_{(I_{1},I_{2},
\ldots,I_{p})} a_{I_{1}}^{j_{1}}a_{I_{2}}^{j_{2}}\ldots a_{I_{p}}^{j_{p}}
y_{j_{1}j_{2}\ldots j_{p}}^{k}},\\
\bar w_{j_{1}j_{2}\ldots j_{s}}^{\sigma} = w_{j_{1}j_{2}\ldots
j_{s}}^{\sigma},\quad 0\le s \le r. \endgathered
$$

Equations
$
w_{j_{1}j_{2}\ldots j_{s}}^{\sigma} = c_{j_{1}j_{2}\ldots j_{s}}^{\sigma},
$
where
$c_{j_{1}j_{2}\ldots j_{s}}^{\sigma} \in \R$ are equations of the orbits of
this action, and the functions 
$y^{i}, w_{j_{1}j_{2}\ldots j_{s}}^{\sigma}$ represent a complete system of
real-valued $L_{n}^{r}$-invariants on
$W$.
Moreover, each of these invariants arises by applying a sequence of the
vector fields
$\Delta_{i}$
to the invariants
$w^{\sigma}$.
\endproclaim

Our aim now will be to express the vector fields $\Delta_{i}$
in terms of the adapted charts
$(W,\Phi)$
(Theorem~2).

\proclaim{Corollary 2}
The vector field
$\Delta_{i}$
has an expression

$$
\gathered
\Delta _{i}(J_{0}^{r}\zeta) = {\partial \over \partial y^{i}} +
\sum\limits_{l=0}^{r-1} 
{\sum\limits_{p_{1}\le p_{2}\le \ldots \le p_{l}} w_{p_{1}p_{2}\ldots
p_{l}i}^{\nu}(J_{0}^{r}\zeta)} 
\left( {\partial \over \partial
w_{p_{1}p_{2}\ldots p_{l}}^{\nu}} \right)_{J_{0}^{r}\zeta}\\ +
\sum\limits_{l=1}^{r-1} 
{\sum\limits_{p_{1}\le p_{2}\le \ldots \le p_{l}}
z_{i}^{s}(J_{0}^{r}\zeta)y_{p_{1}p_{2}
\ldots p_{l}s}^{k}(J_{0}^{r}\zeta)} \left( {\partial \over \partial
y_{p_{1}p_{2}\ldots p_{l}}^{k}} \right)_{J_{0}^{r}\zeta}. \endgathered
\tag{(3.13)}
$$
\endproclaim

\demo{\bf Proof}
We proceed by direct computation, using (2.5) and Theorem 2. \enddemo

Note that at every point of its domain, the vector fields $\Delta_{i}$
(3.12) span an
$n$-dimensional
vector subspace of the tangent space of
$\Imm T_{n}^{r-1}Y$,
determined independently of charts. Indeed, if $(V,\psi)$,
and
$(\bar V,\bar \psi)$,
are two charts, then by (2.6),

$$
\bar \Delta_{i} = \bar z_{i}^{s}\bar d_{s} = \bar z_{i}^{s}d_{s} = \bar
z_{i}^{s}\delta_{s}^{p}d_{p} =
\bar z_{i}^{s}y_{s}^{q}z_{q}^{p}d_{p} = \bar z_{i}^{s}y_{s}^{q}\Delta_{p}.
\tag{(3.14)}
$$

In the following corollary we use these vector fields to derive the
transformation properties of the 
functions $w_{p_{1}p_{2}\ldots p_{k}}^{\nu}$.
Denote
$P = (P_{j}^{i}),$
where

$$
P_{j}^{i} = {\partial \bar y^{i} \over \partial y^{j}} + w^{\nu}_{j}
{\partial \bar y^{i} \over 
\partial w^{\nu}}. $$

Taking
$r = 1$
in (3.14), we get

$$
\gathered
\bar \Delta_{i} = {\partial \over \partial \bar y^{i}} + \bar w^{\nu}
{\partial \over \partial \bar w^{\nu}} = 
\bar z^{s}_{i}y^{q}_{s}\Delta_{q} \\
= \bar z^{s}_{i}y^{q}_{s} \left( {\partial \bar y^{i} \over \partial y^{q}}
+ w^{\nu}_{j} {\partial \bar y^{i} 
\over \partial w^{\nu}}\right) {\partial \over \partial \bar y^{j}} +
\bar z^{s}_{i}y^{q}_{s} \left(
{\partial \bar w^{\nu} \over \partial y^{q}} + w^{\lambda}_{q} {\partial
\bar w^{\nu} \over \partial w^{\lambda}}
\right) {\partial \over \partial \bar w^{\nu}},
\endgathered
$$
from which it follows that
$
\delta_{i}^{j} = \bar z^{s}_{i}y^{q}_{s} P^{j}_{i},\bar w^{\nu}_{i} = \bar
z^{s}_{i}y^{q}_{s} \Delta_{q} 
\bar w^{\nu}. $
The first of these conditions implies that the matrix $P$
is regular, and its inverse,
$P^{-1} = Q = (Q^{i}_{j}),$
satisfies

$$
Q^{i}_{j} = \bar z_{j}^{s}y_{s}^{i}.
$$
From the second condition we derive the following formula 
$\bar w^{\nu}_{i} = Q^{q}_{i}\Delta_{q} \bar w^{\nu}.$ 

\proclaim{Corollary 3}
Let
$(V,\psi),\psi = (y^{A})$,
and
$(\bar V,\bar \psi),\bar \psi =(\bar y^{A})$ be two charts on
$Y$,
such that
$V \cap \bar V\not=\emptyset$.
Consider the associated charts
$(V_{n}^{r},\psi_{n}^{r})$
and
$(\bar V_{n}^{r},\bar \psi_{n}^{r})$
and the charts
$(W,\Phi)$
and
$(\bar W,\bar \Phi)$
on
$\Imm T_{n}^{r}Y$.
Let the transformation equations from
$(V,\psi)$
to
$(\bar V,\bar \psi)$
be written in the form

$$
\bar y^{i}=F^{i}(y^{k},w^{\nu}),\quad
\bar w^{\sigma} = F^{\sigma}(y^{k},w^{\nu}). $$
Then the functions $w_{i_{1}i_{2}\ldots i_{k}i_{k+1}}^{\nu}$ obey the
transformation formulas

$$
\bar w_{i_{1}i_{2}\ldots i_{k}i_{k+1}}^{\nu} = Q_{i_{k+1}}^{p}
\Delta_{p}\bar w_{i_{1}i_{2}\ldots i_{k}}^{\nu}. 
\tag{(3.15)}
$$
\endproclaim

\demo{\bf Proof}
By hypothesis, $\det(y_{i}^{k})\not= 0$,\, hence
$\det(\bar y_{i}^{k})\not= 0$.
Therefore, using (3.7) we get
$
\bar w_{i_{1}i_{2}\ldots i_{k}i_{k+1}}^{\nu} = \bar
z_{i_{k+1}}^{s}\delta^{j}_{s}d_{j} \bar w_{i_{1}i_{2}
\ldots i_{k}}^{\nu} =\bar z_{i_{k+1}}^{s}y_{s}^{p}z_{p}^{j}d_{j} \bar
w_{i_{1}i_{2}\ldots i_{k}}^{\nu} = 
Q_{i_{k+1}}^{p} \Delta_{p}\bar w_{i_{1}i_{2}\ldots i_{k}}^{\nu}. $
\enddemo

A point of
$P^{r}_{n}Y$
containing a regular
$(r,n)$-velocity
$J^{r}_{0}\zeta$
is called an
$(r,n)$-{\it contact element}, or an
$r$-{\it contact element} of an
$n$-dimensional
submanifold of
$Y$,
and is denoted by
$[J^{r}_{0}\zeta].$
As in the case of
$r$-jets, the point
$0 \in \R^{n}$
(resp.
$\zeta(0) \in Y$)
is called the {\it source} (resp. the {\it target}) of $[J^{r}_{0}\zeta].$
The set
$G^{r}_{n}$
of
$(r,n)$-contact
elements with source
$0 \in \R^{n}$
and target
$0 \in \R^{n+m}$,
endowed with the natural smooth structure, is called the $(r,n)$-{\it
Grassmannian},
or simply a {\it higher order Grassmannian}. It is standard to check that
the manifold
$P^{r}_{n}Y = \Imm T^{r}_{n}Y/L^{r}_{n}$ is a fiber bundle over
$Y$
with fiber
$G^{r}_{n}$.
$P^{r}_{n}Y$
with this structure is called the
$(r,n)$-{\it Grassmannian bundle},
or simply a {\it higher order Grassmannian bundle} over
$Y$.

Besides the quotient projection
$\rho^{r}_{n}: \Imm T^{r}_{n}Y \to P^{r}_{n}$ (Corollary 1 to Theorem 1) we
have for every $s, 0 \le s \le r$,
the {\it canonical projection} of
$P^{r}_{n}Y$
onto
$P^{s}_{n}Y$
defined by
$\rho^{r,s}_{n}([J^{r}_{0}\zeta]) = [J^{s}_{0}\zeta].$ 

Now we are going to introduce some charts on the manifold of contact
elements $P^{r}_{n}Y$.
To this purpose we consider the {\it adapted charts} on $
\Imm T^{r}_{n}Y,(W,\Phi)$,

$$\Phi = (w^{\sigma},w^{\sigma}_{p_{1}},w^{\sigma}_{p_{1}p_{2}},\ldots,
w^{\sigma}_{p_{1}p_{2}\ldots p_{r}},
y^{i},y^{i}_{j_{1}},y^{i}_{j_{1}j_{2}}, \ldots, y^{i}_{j_{1}j_{2}\ldots
j_{r}}), $$
introduced in Theorem 2. We denote
$\tilde W = \rho^{r}_{n}(W)$,
and if
$J^{r}_{0}\zeta \in W$,
we define
$$
\tilde \Phi = (\tilde y^{i},\tilde w^{\sigma},\tilde w^{\sigma}_{j_{1}},
\tilde w^{\sigma}_{j_{1}j_{2}},\ldots,
\tilde w^{\sigma}_{j_{1}j_{2}\ldots j_{r}}) $$
by

$$
\tilde y^{i}([J^{r}_{0}\zeta]) = y^{i}(J^{r}_{0}\zeta), \quad \tilde
w^{\sigma}_{j_{1}j_{2}\ldots j_{k}}
([J^{r}_{0}\zeta]) = w^{\sigma}_{j_{1}j_{2}\ldots j_{k}}(J^{r}_{0}\zeta).
\tag{(3.16)}
$$

Then the pair
$(\tilde{W},\tilde{\Phi})$
is the {\it associated chart} on
$P^{r}_{n}Y$.
In terms of
$(W,\Phi)$
and
$(\tilde{W},\tilde{\Phi})$
the quotient projection
$\rho^{r}_{n}$
is expressed by the equations

$$
\tilde y^{i} \circ \rho^{r}_{n} = y^{i},\quad
\tilde{w}^{\sigma}_{j_{1}j_{2}\ldots j_{k}} \circ \rho^{r}_{n} = 
w^{\sigma}_{j_{1}j_{2}\ldots j_{k}}.
\tag{(3.17)}
$$

Consider a point
$J^{r}_{0}\zeta \in W$,
and the vector subspace of the tangent space
$T_{\rho^{r}_{n}(J^{r}_{0}\zeta)}P^{r}_{n}Y$ spanned by the vectors
$T_{J^{r}_{0}\zeta}\rho^{r}_{n}\cdot \Delta_{i}(J^{r}_{0}\zeta),$ where the
vectors
$\Delta_{i}(J^{r}_{0}\zeta)$
are defined by (3.13) and (2.5). Indeed, this vector subspace is
independent of the choice of a chart used in the 
definition of $d_{i}.$
It follows from (3.13) and (3.17) that the vector field $\Delta_{i}$
is
$\rho^{r}_{n}$-projectable, and its
$\rho^{r}_{n}$-projection is the vector field 

$$
\tilde{\Delta}_{i} = {\partial \over \partial \tilde y^{i}} +
\sum\limits_{p=0}^{r-1} {\sum\limits_{j_{1}
\le j_{2}\le \ldots \le j_{p}} \tilde w_{j_{1}j_{2}\ldots j_{p}i}^{\sigma}}
{\partial \over \partial
\tilde w_{j_{1}j_{2}\ldots j_{p}}^{\sigma}}. $$

Thus we have the following commutative diagram: 

$$
\CD
\Imm T^{r}_{n} @>\Delta_{i}>> T \Imm T^{r-1}_{n} \\ @VV\rho^{r}_{n}V	@VVT
\rho^{r}_{n}V \\
P^{r}_{n}Y @>\tilde\Delta_{i}>> T P^{r-1}_{n}Y \endCD
$$

From now on we adopt the standard convention for writing fibered
coordinates, 
and we omit the tilde over the coordinate functions on the left in (3.16).
Then the coordinate functions of the 
chart
$(\tilde{W},\tilde{\Phi})$
will be denoted simply by
$
\tilde \Phi = (y^{i},w^{\sigma},w^{\sigma}_{j_{1}},w^{\sigma}_{j_{1}j_{2}},
\ldots,w^{\sigma}_{j_{1}j_{2}
\ldots j_{r}}). $

Let us consider two charts
$(V,\phi), \phi = (y^{A})$,
and
$(\bar V,\bar \phi), \bar \phi = (\bar y^{A})$, such that
$V \cap \bar V \not= 0,$
and the associated charts
$(V^{r}_{n},\phi^{r}_{n})$
and
$(\bar V^{r}_{n},\bar \phi^{r}_{n})$
on
$\Imm T^{r}_{n}Y.$
The transformation equations for the corresponding associated charts on
$P^{r}_{n}Y$
are given by
$
\bar w^{\nu}_{i_{1}i_{2}\ldots i_{k}i_{k+1}} = Q^{p}_{i_{k+1}}\Delta_{p}
\bar w^{\nu}_{i_{1}i_{2}\ldots i_{k}}
$
(3.15).

\heading 4. Scalar invariants of $(r,n)$-velocities \endheading 

Our aim in this section will be to describe all continuous
$L_{n}^{r}$-invariant, real-valued functions on the 
manifold of $(r,n)$-velocities
$T_{n}^{r}Y$.

As in the case of regular
$(r,n)$-velocities, we denote by
$\rho^{r}_{n}: T^{r}_{n}Y \to T^{r}_{n}Y /L^{r}_{n}$ the {\it canonical
quotient projection}. The quotient set 
$T^{r}_{n}Y /L^{r}_{n}$
will be considered with its {\it canonical topological structure}; then
$\rho^{r}_{n}$
is an {\it open} mapping. The set
$\Imm T^{r}_{n}$
is an open, dense, subset of
$T^{r}_{n}.$
We have the canonical projection
$\pi^{r}_{n}: T^{r}_{n}Y/L^{r}_{n} \to Y$, as well as its restriction
$\pi^{r}_{n}: \Imm T^{r}_{n}Y /L^{r}_{n} \to Y$ to the
$(r,n)$-{\it Grassmann bundle}
$P^{r}_{n} = \Imm T^{r}_{n}Y/L^{r}_{n}$, which are both continuous. These
mappings define a commutative diagram 

$$
\CD
\Imm T_{n}^{r}Y @> >> T_{n}^{r-1}Y \\
@V VV @V VV \\
P_{n}^{r}Y @> >> T_{n}^{r}Y/L^{r}_{n}\\
@V VV @V VV \\
Y @> >> Y
\endCD
$$

$P_{n}^{r}Y$
is an open, dense subset of
$T_{n}^{r}Y/L^{r}_{n}.$
Indeed
$P_{n}^{r}Y$
is open in
$T_{n}^{r}Y/L^{r}_{n}$
by the definition of the quotient topology, since $\Imm T_{n}^{r}Y =
(\rho^{r}_{n})^{-1}(P^{r}_{n}Y)$ is open in
$T_{n}^{r}Y.$
If
$[J^{r}_{0}\chi_{0}] \in T^{r}_{n}Y/L^{r}_{n}$ is such that
$[J^{r}_{0}\chi_{0}] \not\in P^{r}_{n}Y$, and
$W$
is a neighbourhood of
$[J^{r}_{0}\chi_{0}]$,
then
$(\rho^{r}_{n})^{-1}(W)$
is an open set in
$T^{r}_{n}Y$
containing
$[J^{r}_{0}\chi_{0}]$
as a subset. Since
$\Imm T_{n}^{r}Y$
is dense in
$T^{r}_{n}Y$,
$(\rho^{r}_{n})^{-1}(W)\cap \Imm T^{r}_{n}Y$ is a non-empty open subset of
$\Imm T_{n}^{r}Y$,
and since
$\rho^{r}_{n}$
is open, the set
$\rho^{r}_{n}((\rho^{r}_{n})^{-1}(W)\cap \Imm T^{r}_{n}Y)$ is open in
$P_{n}^{r}Y.$
But
$\rho^{r}_{n}((\rho^{r}_{n})^{-1}(W)\cap \Imm T^{r}_{n}Y) \subset W$ which
means that
$W$
contains an element of the set
$P_{n}^{r}Y.$

Any continuous function on a subset of
$P_{n}^{r}Y$
defines, when composed with the quotient projection $\rho^{r}_{n}: \Imm
T^{r}_{n}Y \to P_{n}^{r}Y$, an
$L_{n}^{r}$-invariant, continuous function on the corresponding subset of
$\Imm T^{r}_{n}Y$,
and vice versa, any
$L_{n}^{r}$-invariant, continuous function on an open,
$L_{n}^{r}$-invariant subset of
$P_{n}^{r}Y$
can be factored through
$\rho^{r}_{n}.$ Since the values of a continuous, real-valued function on 
$T_{n}^{r}Y/L_{n}^{r}$
are uniquely determined by its values on $P_{n}^{r}Y,$
the projection
$\rho^{r}_{n}$
is the {\it basis of}
$L_{n}^{r}$-{\it invariant functions} on $T_{n}^{r}Y.$

It is now clear that our problem of finding all continuous
$L_{n}^{r}$-invariant, real-valued function on 
$T_{n}^{r}Y$
is equivalent with the problem of finding continuous functions on open
subset of the quotient
$T_{n}^{r}Y/L_{n}^{r}.$
This gives rise to the problem of {\it continuous prolongation} of
functions on $P_{n}^{r}Y$
to the quotient space
$T_{n}^{r}Y/L_{n}^{r}.$

First we need to discuss separability of the points on
$T_{n}^{r}Y/L_{n}^{r}.$
It is easily seen that the quotient topology on $T_{n}^{r}Y/L_{n}^{r}$
is {\it not} Hausdorff.

Note that any two points
$[J_{0}^{r}\chi_{0}] ,[J_{0}^{r}\chi] \in T_{n}^{r}Y/L_{n}^{r}Y$ such that
$\chi_{0}(0) \not= \chi(0)$,
can always be separated by open sets. This follows from the continuity of
the quotient projection of
$\pi^{r}_{n}$,
and from separability of
$Y.$
To study the situation in the fibers, we prove the following lemma. 

\proclaim{Lemma 4}
Let
$y \in Y$
be a point,
$(V,\psi), \psi = (y^{A})$
a chart at
$y$,
and
$J^{r}_{0}\chi_{0} \in (\tau^{r,0}_{n})^{-1}(y)$ the
$(r,n)$-velocity with target
$y$
defined by
$J^{r}_{0}\chi_{0} = (y^{A},0,0,\ldots,0)$ in the associated chart. Then
any
$L_{n}^{r}$-invariant neighbourhood of
$J^{r}_{0}\chi_{0}$ contains the fiber
$(\tau^{r.0}_{n})^{-1}(y).$
\endproclaim

\demo{\bf Proof}
The fiber
$(\tau^{r,0}_{n})^{-1}(y) =
(\rho^{r,0}_{n})^{-1}((\tau^{r,0}_{n})^{-1}(y))$ over
$y \in Y$
in
$T_{n}^{r}Y$
is endowed with the induced chart
$
(V^{r}_{n},\psi^{r}_{n}), 
\psi^{r}_{n} = (y^{A},y^{A}_{j_{1}},y^{A}_{j_{1}j_{2}},\ldots,
y^{A}_{j_{1}j_{2}\ldots j_{r}}).
$
The coordinates
$
\bar y^{A}, \bar y^{A}_{j_{1}j_{2}\ldots j_{s}}$ of the points of the orbit
$[J^{r}_{0}\chi_{0}]$
are given by (2.16),

$$
\bar y^{A} = y^{A}, \quad
\bar y^{A}_{j_{1}j_{2}\ldots j_{s}} =
\sum\limits_{p=1}^{s} {\sum\limits_{(I_{1},I_{2},\ldots I_{p})}
a^{j_{1}}_{I_{1}}a^{j_{2}}_{I_{2}}
\ldots a^{j_{p}}_{I_{p}} y^{A}_{j_{1}j_{2}\ldots j_{p}}},
$$
where
$
J^{r}_{0}\alpha \in L^{r}_{n},
J^{r}_{0}\alpha = (a^{i}_{j_{1}},a^{i}_{j_{1}j_{2}},\ldots,
a^{i}_{j_{1}j_{2}\ldots j_{r}}).
$
Thus,
$\bar y^{A} = y^{A},  \bar y^{A}_{j_{1}j_{2}\ldots j_{s}} = 0,$ which means
that the orbit
$[J^{r}_{0}\chi_{0}]$
consists of a single point.
Let
$W$
be an
$L^{r}_{n}$-invariant
neighbourhood of the point
$J^{r}_{0}\chi_{0}.$
We show that each orbit in
$(\tau^{r,0}_{n})^{-1}(y)$
has a non-empty intersection with
$W.$
Then we apply
$L^{r}_{n}$-invariance to obtain the inclusion $(\tau^{r,0}_{n})^{-1}(y)
\subset W.$

Let us consider an arbitrary element
$
J^{r}_{0}\chi = (y^{A},y^{A}_{j_{1}},y^{A}_{j_{1}j_{2}},\ldots,
y^{A}_{j_{1}j_{2}\ldots j_{r}}) 
\in (\tau^{r,0}_{n})^{-1}(y), $
and a one-parameter family of velocities $
J^{r}_{0}\chi \circ J^{r}_{0}\beta_{\tau} $
in
$(\tau^{r,0}_{n})^{-1}(y)$
defined in components by 
$$
\beta_{\tau} = (\beta_{\tau}^{i}),
\beta_{\tau}^{i}(t^{1},t^{2},\ldots,t^{n}) = \tau t^{i}, $$
where
$0 \le \tau \le 1.$
Then
$J^{r}_{0}\beta_{\tau} = (\tau \delta^{i}_{j},0,0,\ldots,0),$ and by
(2.16), the
$L^{r}_{n}$-orbit of
$J^{r}_{0}\chi$
contains the points
$
J^{r}_{0}\chi \circ J^{r}_{0}\beta_{\tau} = (\bar y^{A},\bar y^{A}_{j_{1}},
\bar y^{A}_{j_{1}j_{2}},\ldots,
\bar y^{A}_{j_{1}j_{2}\ldots j_{r}}) $
given by
$
\bar y^{A} = y^{A}, \bar y^{A}_{j_{1}j_{2}\ldots j_{s}} = \tau^{s} \bar
y^{A}_{j_{1}j_{2}\ldots j_{s}}. $
Clearly, for sufficiently small
$\tau, J^{r}_{0}\chi \circ J^{r}_{0}\beta_{\tau} \in W.$ 

This shows that the orbits passing through any neighbourhood of the point
$J^{r}_{0}\chi_{0}$,
where
$\chi_{0}(0) = y$,
fill the whole fiber
$(\tau^{r,0}_{n})^{-1}(y).$
\enddemo

Consider a point
$y \in Y$,
a chart
$(V,\psi), \psi = (y^{A})$
at
$y$,
and the
$L^{r}_{n}$-orbit
$[J^{r}_{0}\chi_{0}]$
of the velocity
$J^{r}_{0}\chi_{0} = (y^{A},0,0,\ldots,0) \in (\tau^{r,0}_{n})^{-1}(y).$
Lemma 4 shows that any neighbourhood 
of the orbit $[J^{r}_{0}\chi_{0}] \in T^{r}_{n}Y/L^{r}_{n}$ contains the
fiber
$(\tau^{r,0}_{n})^{-1}(y)$
in
$T^{r}_{n}Y/L^{r}_{n}$
over
$y.$
This proves, in particular, that no point of $(\tau^{r}_{n})^{-1}(y)$
can be separated from
$[J^{r}_{0}\chi_{0}]$
by open sets.

This gives us the following theorem saying that if a continuous invariant
is defined on a fiber in
$T^{r}_{n}Y$,
then it is constant along this fiber.

\proclaim{Theorem 3}
Let
$W$
be an open,
$L_{n}^{r}$-invariant set in
$\Imm T_{n}^{r}Y, f:W\to \R$
an
$L_{n}^{r}$-invariant function. Assume that $W$
contains two regular velocities
$J^{r}_{0}\zeta, J^{r}_{0}\chi$
with common target
$y = \zeta(0) = \chi(0),$
such that
$f(J^{r}_{0}\zeta) \not= f(J^{r}_{0}\chi).$ Then
$f$
cannot be continuously prolonged to the fiber $\tau_{n}^{r,0}(y) \subset
\Imm T_{n}^{r}Y$. \endproclaim

\demo{Proof}
Indeed, since
$\R$
is Hausdorff, any continuous,
$L_{n}^{r}$-invariant,
real-valued function takes the same value at the points which cannot be
separated by open sets. Assume that
$f$
can be prolonged to the fiber
$\tau_{n}^{r,0}(y) \subset \Imm T_{n}^{r}Y.$ Then by lemma 4,
$f$
is equal along the fiber
$\tau_{n}^{r,0}(y)$
to
$f(J^{r}_{0}\chi_{0}) = const$,
which is a contradiction.
\enddemo

In particular, none of the
$L_{n}^{r}$-invariant function
$w^{\sigma}_{j_{1}j_{2}\ldots j_{k}}$
(Theorem 2) can be prolonged to a fiber
$(\tau_{n}^{r,0})^{-1}(y).$

Now it is immediate that each
$L_{n}^{r}$-invariant function on
$T^{r}_{n}Y$
is trivial in the following sense.

\proclaim{Corollary 1}
A, continuous function
$f:T_{n}^{r}Y \to \R$
is
$L_{n}^{r}$-invariant if and only if
$f = F \circ \tau_{n}^{r,0},$
where
$F:Y \to \R$
is a continuous function.
\endproclaim

\heading 5. Appendix: Regular $(2,n)$-velocities \endheading 

As before,
$Y$
denotes a smooth manifold of dimension
$n + m.$
In this section we consider the manifold $\Imm T^{2}_{n}Y$
of regular
$(2,n)$-velocities on
$Y$,
and the Grassmann bundle
$P^{2}_{n}Y.$
We wish to collect in an explicit form all basic formulas concerning charts
and invariants in this case, 
which will be important for applications. 

If
$(V,\psi), \psi = (y^{A}),$
is a chart on
$Y$,
define
$V^{2}_{n} = (\tau^{2,0}_{n})^{-1}(V)$
and
$\psi^{2}_{n} = (y^{A},y^{A}_{i},y^{A}_{ij})$ where
$1 \le A \le n+m, 1 \le i \le j \le n$,
by the formulas
$
y^{A}(J^{2}_{0}\zeta) = y^{A}(\zeta(0)), y^{A}_{i}(J^{2}_{0}\zeta) =
D_{i}y^{A}(\zeta(0)), 
y^{A}_{ij}(J^{2}_{0}\zeta) = D_{i}D_{j}y^{A}(\zeta(0)).$ If
$(\bar V,\bar \psi), \bar \psi = (\bar y^{A}),$ is another chart on
$Y$,
and the transformation equations are written as $\bar y^{A} =
F^{A}(y^{B})$,
then

$$
\bar y^{A} = F^{A}(y^{B}), \quad
\bar y^{A}_{i} = {\partial F^{A} \over \partial y^{B}} y^{B}_{i}, \quad
\bar y^{A}_{ij} = 
{\partial^{2} F^{A}\over \partial Y^{B} \partial y^{C}} y^{B}_{i}y^{C}_{j}
+ {\partial F^{A} \over 
\partial y^{B}} y^{B}_{ij} \tag{(5.1)}
$$
on
$V^{2}_{n} \cap \bar V^{2}_{n}$
(see (2.9)-(2.11)).

If
$f: V^{1}_{n} \to \R$
is a smooth function, we define a function $d_{i}f: V^{2}_{n} \to \R$
by

$$
d_{i}f = {\partial f \over \partial y^{A}} y^{A}_{i} + {\partial f \over
\partial y^{A}_{j}} y^{A}_{ij} $$
(the
$i$-th
{\it formal derivative} of
$f$).
In particular,
$d_{i}y^{A} = y^{A}_{i}, d_{i}y^{A}_{j} = y^{A}_{ij}.$ 

By definition,
$\text{rank}(y^{B}_{s}(J^{2}_{n}\zeta)) = n$ at every point
$J^{2}_{n}\zeta \in V^{2}_{n}.$
Thus, there exists a subsequence
$I = (A_{1},A_{2},\ldots,A_{n})$
of the sequence
$(1,2,\ldots,n,n+1,n+m)$
such that
$\det(y^{A_{i}}_{j}(J^{2}_{n}\zeta)) \not= 0.$ Denote
$V^{2(I)}_{n} = \{ J^{2}_{0}\zeta \in V^{2}_{n} |
\det(y^{A_{i}}_{j}(J^{2}_{n}\zeta)) \not= 0 \}.$ If
$\psi^{2(I)}_{n}$
is the restriction of
$\psi^{2}_{n}$
to
$V^{2(I)}_{n}$,
then the pair
$
(V^{2(I)}_{n},\psi^{2}_{n}), \psi^{2(I)}_{n} =
(y^{A},y^{A}_{i},y^{A}_{ij}),
$
is a chart on
$\Imm T^{2}_{n}Y$,
and

$$
\bigcup_{I} V^{2(I)}_{n} = V^{2}_{n}.
$$

By (2.14), the group multiplication
$(J^{2}_{n}\alpha,J^{2}_{n}\beta) \to J^{2}_{n}\alpha \circ J^{2}_{n}\beta$
in the second differential group
$L^{2}_{n}$
of
$\R^{n}$
is given in the {\it canonical coordinates} by 

$$
c^{k}_{i} = b^{p}_{i}a^{k}_{p}, \quad
c^{k}_{ij} = b^{p}_{ij}a^{k}_{p} + b^{p}_{i}b^{q}_{j}a^{k}_{pq}.
\tag{(5.2)}
$$
Indeed, in this formula
$a^{k}_{p}, a^{k}_{pq}$
(resp.
$b^{k}_{p}, b^{k}_{pq}$,
resp.
$c^{k}_{p}, c^{k}_{pq}$)
are the coordinates of
$J^{2}_{0}\alpha$
(resp.
$J^{2}_{0}\beta$,
resp.
$J^{2}_{0}\alpha \circ J^{2}_{0}\beta$). $L^{2}_{n}$
acts on
$\Imm T^{2}_{n}Y$
smoothly to the right by the jet composition $
(J^{2}_{0}\zeta \circ J^{2}_{0}\alpha) \to J^{2}_{0}\zeta \circ
J^{2}_{0}\alpha.
$
This action is expressed by

$$
\bar y^{A} = y^{A}, \quad
\bar y^{A}_{i} = y^{A}_{s}a^{s}_{i}, \quad \bar y^{A}_{ij} =
y^{A}_{pq}a^{p}_{i}a^{q}_{j} + y^{A}_{p}a^{p}_{ij},
\tag{(5.3)}
$$
where
$y^{A}, y^{A}_{i}, y^{A}_{ij}$
are the coordinates of a velocity
$J^{2}_{0}\zeta$,
and
$\bar y^{A}, \bar y^{A}_{i},\bar y^{A}_{ij}$ are the coordinates of the
transformed velocity $J^{2}_{0}\zeta 
\circ J^{2}_{0}\alpha.$

Now we are going to construct an atlas on $\Imm T^{2}_{n}Y$,
adapted to this group action. Given a chart $(V,\psi), \psi = (y^{A}),$
we note that the action (5.3) preserves each of the sets $V^{2(I)}_{n}.$
Indeed, if
$J^{2}_{0}\zeta \in V^{2(I)}_{n}$,
then by definition, the matrix
$y^{A_{i}}_{j} = y^{A_{i}}_{j}(J^{2}_{0}\zeta)$ is of maximal rank, and the
second equation (5.3) implies that 
the matrix $\bar y^{A_{i}}_{j} = y^{A_{i}}_{j}(J^{2}_{0}\zeta\alpha)$ of
the transformed point is also of maximal 
rank. 

Consider for example the case
$I = (1,2,\ldots ,n).$
Then
$\det(y^{i}_{j}) \not= 0$
for
$1 \le i,j \le n$
(i.e. on
$V^{2(I)}_{n}$).
We define smooth functions
$z^{i}_{j}: V^{2(I)}_{n} \to \R$
by
$z^{i}_{p}y^{p}_{j} = \delta^{i}_{j}.$
These functions form a regular matrix on $V^{2(I)}_{n}.$
Equations (5.3) then give for
$A = k = 1,2,\ldots ,n$

$$
\gathered
z^{k}_{p}\bar y^{p}_{i} = a^{k}_{i},\quad
z^{k}_{s}\bar y^{s}_{ij} = z^{k}_{s}y^{s}_{pq}a^{p}_{i}a^{q}_{j} +
z^{k}_{s}y^{s}_{p}a^{p}_{ij} = 
z^{k}_{s}y^{s}_{pq}z^{p}_{r} \bar y^{r}_{i} z^{q}_{t} \bar y^{t}_{j} +
a^{k}_{ij}, \endgathered
$$
and for
$A = \sigma = n+1,n+2,\ldots , m$

$$
\gathered
\bar y^{A} = y^{A}, \quad \bar y^{A}_{i} = y^{A}_{s}z^{s}_{p}\bar
y^{p}_{i},\quad \bar y^{A}_{ij} = 
y^{A}_{pq}z^{p}_{s} \bar y^{s}_{i} z^{q}_{t} \bar y^{t}_{j} +
y^{A}_{k}(z^{k}_{s} \bar y^{s}_{ij}-z^{k}_{s} 
\bar y^{s}_{pq}z^{p}_{r}\bar y^{r}_{i} z^{q}_{t} \bar y^{t}_{j}).
\endgathered
$$
Since the second formula gives us the relation $\bar z^{i}_{j} \bar
y^{\sigma}_{i} = 
z^{s}_{j} y^{\sigma}_{s}$, and the third one implies

$$
\bar z^{i}_{u}\bar z^{j}_{v}(\bar y^{\sigma}_{ij}-\bar z^{k}_{s} \bar
y^{\sigma}_{k}\bar y^{s}_{ij}) = 
y^{\sigma}_{pq}z^{p}_{u} z^{q}_{v}- z^{k}_{s}
y^{\sigma}_{k}y^{s}_{pq}z^{p}_{u}z^{q}_{v} = 
z^{p}_{u} z^{q}_{v}(y^{\sigma}_{pq}-z^{k}_{s} y^{\sigma}_{k}y^{s}_{pq}),
$$
we finally get

$$
\bar y^{A}= y^{A},\quad
\bar z^{i}_{j} \bar y^{\sigma}_{i} = z^{s}_{j} y^{\sigma}_{s}, $
and
$
\bar z^{i}_{u}\bar z^{j}_{v}(\bar y^{\sigma}_{ij}-\bar z^{k}_{s} \bar
y^{\sigma}_{pq}\bar y^{s}_{ij}) = 
z^{p}_{u} z^{q}_{v} (y^{\sigma}_{pq}-z^{k}_{s}y^{\sigma}_{k}y^{s}_{pq}). $$
Therefore, the functions

$$
y^{i}, \quad w^{\sigma} = y^{\sigma}, \quad w^{\sigma}_{i} =
z^{s}_{i}y^{\sigma}_{s}, \quad w^{\sigma}_{ij} = 
z^{p}_{i} z^{q}_{j}
(y^{\sigma}_{pq}-z^{k}_{s}y^{\sigma}_{k}y^{s}_{pq}) \tag{(5.4)}
$$
are {\it constant} along the
$L^{2}_{n}$-orbits in
$V^{2(I)}_{n} \subset \Imm T^{2}_{n}Y$,
and the functions
$y^{i}, w^{\sigma}, w^{\sigma}_{ij}, y^{i}_{j}, y^{i}_{jk}$ define a new
chart on the set
$V^{2(I)}_{n}$
{\it adapted} to the group action of
$L^{2}_{n}$.
The right action (5.3) is expressed in terms of this new chart by $
\bar y^{i} = y^{i}, \bar w^{\sigma} = w^{\sigma}, \bar w^{\sigma}_{ij} =
w^{\sigma}_{ij}, \bar y^{k}_{i} = 
y^{k}_{s}a^{s}_{i}, \bar y^{k}_{ij} = y^{k}_{pq}a^{p}_{i}
a^{q}_{j} + y^{k}_{p}a^{p}_{ij}.
$
The functions (5.4) are components of the quotient projection
$\rho^{2}_{n}$,
of the manifold of regular
$(2,n)$-velocities
$\Imm T^{2}_{n}Y$
onto the
$(2,n)$-Grassmann bundle
$P^{2}_{n}Y$,
and form a {\it basis of}
$L^{2}_{n}$-{\it invariant functions} on $V^{2(I)}_{n}.$

A direct interpretation of the coordinate functions (5.4) is obtained as
follows. Let
$(V,\psi), \psi = (x^{i},y^{\sigma})$,
be a chart on
$Y$,
and let
$J^{2}_{0}\zeta \in V^{2(I)}_{n}.$
We assign to the
$2$-jet
$J^{2}_{0}\zeta$
an element
$J^{2}_{0}\alpha \in L^{2}_{n}$
by means of a representative
$\alpha$
satisfying, in addition to the condition $\alpha(0) = 0$,
the following two conditions
$
a^{s}_{i}(J^{2}_{0}\alpha) = D_{i}\alpha^{s}(0) =
y^{s}_{i}(J^{2}_{0}\zeta), a^{s}_{ij}(J^{2}_{0}\alpha) = 
D_{i}D_{j}\alpha^{s}(0) = y^{s}_{ij}(J^{2}_{0}\zeta).
$
Then it is easily seen that

$$
\gathered
a^{i}_{r}(J^{2}_{0}\alpha^{-1}) = D_{r}(\alpha^{-1})^{i}(0) =
z^{i}_{r}(J^{2}_{0}\zeta),\\
a^{i}_{rs}(J^{2}_{0}\alpha^{-1}) = D_{r}D_{s}(\alpha^{-1})^{i}(0) = 
- z^{j}_{r}(J^{2}_{0}\zeta) z^{k}_{s}a^{i}_{r}(J^{2}_{0}\zeta)
y^{p}_{jk}(J^{2}_{0}\zeta),
\endgathered
$$
and we get for the coordinates of the
$2$-jet
$J^{2}_{0}(\zeta\circ \alpha^{-1})\in V^{2(I)}_{n}$, by (5.4),

$$
\gathered
y^{i}(J^{2}_{0}(\zeta\circ\alpha^{-1})) = y^{i}(J^{2}_{0}\zeta), \quad
y^{k}_{i}(J^{2}_{0}
(\zeta\circ\alpha^{-1})) = \delta^{k}_{i}, \quad
y^{k}_{ij}(J^{2}_{0}(\zeta\circ\alpha^{-1})) = 
0,\\ y^{\sigma}(J^{2}_{0}(\zeta\circ\alpha^{-1})) =
w^{\sigma}(J^{2}_{0}\zeta), \quad y^{\sigma}_{i}(J^{2}_{0}
(\zeta\circ\alpha^{-1})) = w^{\sigma}_{i}(J^{2}_{0}\zeta),\\
y^{\sigma}_{ij}(J^{2}_{0}(\zeta\circ\alpha^{-1})) =
w^{\sigma}_{ij}(J^{2}_{0}\zeta),
\endgathered
$$
where
$k = 1,2,\ldots,n \quad \sigma = n+1,n+2,\ldots ,m.$ This represents the
desired interpretation of the functions 
(5.4) as jet coordinates of the
$2$-jets
$J^{2}_{0}\zeta \circ J^{2}_{0}\alpha$,
with
$J^{2}_{0}\alpha$
determined by the considered chart.

One can determine the transformation formulas from
$\psi^{2}_{n}(J^{2}_{0}(\zeta\circ\alpha^{-1}))$ to
$\bar \psi^{2}_{n}(J^{2}_{0}(\zeta\circ\bar\alpha^{-1}))$, with obvious
meaning of
$\bar\alpha$.
These formulas illustrate the well-known fact that the higher order
Grassmann bundles have a relatively 
complicated smooth structure. Consider an element $J^{2}_{0}\zeta \in
V^{2(I)}_{n} \cap \bar V^{2(I)}_{n}.$ 
The corresponding computations for
$J^{2}_{0}\zeta \in V^{2(J)}_{n} \cap \bar V^{2(J)}_{n}$ with arbitrary
$I, J$
are quite analogous. We have

$$
\gathered
\bar \psi^{2}_{n}(J^{2}_{0}(\zeta\circ\bar\alpha^{-1})) = \bar
\psi^{2}_{n}(J^{2}_{0}(\zeta\circ\alpha^{-1}\circ 
\alpha\bar\alpha^{-1})) = \bar
\psi^{2}_{n}(J^{2}_{0}(\zeta\circ\alpha^{-1}))\circ J^{2}_{0}
(\alpha\bar\alpha^{-1}))\\
= \bar \psi^{2}_{n}(\psi^{2}_{n})^{-1}
(\psi^{2}_{n}(J^{2}_{0}(\zeta\circ\alpha^{-1})\circ
J^{2}_{0}(\alpha\bar\alpha^{-1}))).
\endgathered
$$
To derive explicit expressions, one has to substitute for the group
multiplication (5.2), the group action (5.3), 
and the transformation (5.1) in this formula.
Denoting

$$
P^{q}_{s} = {\partial F^{q}\over \partial y^{s}} + {\partial F^{q}\over
\partial y^{\nu}}w^{\nu}_{s}, $$
we get a regular matrix
$P = (P^{q}_{s}).$
Let
$Q = P^{-1} = (Q^{i}_{j})$
be its inverse. Then

$$
\bar y^{q}_{i} = P^{q}_{s}y^{s}_{i},\quad \bar z^{q}_{i} =
Q^{s}_{i}z^{q}_{s}. $$
After a tedious but straightforward calculation we get the following
result. Given the transformation equations 
on
$
Y, \bar y^{k} = F^{k}(y^{p},y^{\nu}),\bar y^{\sigma} =
F^{\sigma}(y^{p},y^{\nu}),
$
then on
$V^{2(I)}_{n} \cap \bar V^{2(J)}_{n},$

$$
\gathered
\bar y^{k} = F^{k}(y^{p},y^{\nu}),\quad
\bar w^{\sigma} = F^{\sigma}(y^{p},y^{\nu}),\quad \bar w^{\sigma}_{i} =
Q^{p}_{i} \left(
{\partial F^{\sigma}\over \partial y^{p}} + {\partial F^{\sigma}\over
\partial w^{\nu}} w^{\nu}_{p} \right), \\
\bar w^{\sigma}_{ij} = Q^{p}_{i}Q^{q}_{j} \left( {\partial^{2}
F^{\sigma}\over \partial y^{p}\partial y^{q}} + 
{\partial^{2} F^{\sigma}\over \partial y^{p}\partial w^{\nu}} w^{\nu}_{q}+
{\partial^{2} F^{\sigma}\over \partial
 w^{\nu}\partial y^{q}} w^{\nu}_{p} + {\partial^{2} F^{\sigma}\over
\partial w^{\mu}\partial w^{\nu}} 
w^{\mu}_{p}w^{\nu}_{q}
+ {\partial F^{\sigma}\over \partial w^{\nu}} w^{\nu}_{pq}\right) -
Q^{a}_{i}Q^{b}_{j}Q^{p}_{k} \\
\left( {\partial F^{\sigma}\over \partial y^{p}} + {\partial
F^{\sigma}\over \partial w^{\nu}} w^{\nu}_{p}\right) 
\left( {\partial^{2} F^{k}\over \partial y^{a}\partial y^{b}} +
{\partial^{2} F^{k}\over \partial y^{a}\partial 
w^{\sigma}} w^{\sigma}_{b}+ {\partial^{2} F^{k}\over \partial
w^{\sigma}\partial y^{b}} w^{\sigma}_{a} + 
{\partial^{2} F^{k}\over \partial w^{\sigma}\partial w^{\lambda}}
w^{\sigma}_{a}w^{\lambda}_{b} +
{\partial F^{k}\over \partial w^{\lambda}} w^{\lambda}_{ab} \right).
\endgathered
$$

Clearly, these equations represent the transformation formulas for the
induced charts on the
$(2,n)$-Grassmann bundle
$P^{2}_{n}Y.$

\heading References \endheading
\leftskip 1.5pc

\def\bibitem[#1]{\noindent\hglue-1pc\hbox to 1pc{\hss[#1]\ }\ignorespaces}

\bibitem[1] P. Dedecker, On the generalization of symplectic geometry to
multiple integrals in the calculus 
of variations, Lecture Notes in Math. 570, Springer, Berlin, 395--456.

\bibitem[2] J. Dieudonn\'e, {\it Elements d'Analyse 3}, Gauthier-Villars,
Paris, 1970. 

\bibitem[3] D.R. Grigore, A generalized Lagrangian formalism in particle
mechanics and classical field theory, 
Fortschr. d. Phys. 41 (1993), 569--617. 

\bibitem[4] D.R. Grigore and O. Popp, On the Lagrange-Souriau form in
classical field theory, submitted 
for publication

\bibitem[5] P.. Horv\'athy, Variational formalism for spinning particles,
JMP 20 (1974), 49--52.

\bibitem[6] J. Klein, Espaces variationels et Mechanique, Ann. Inst.
Fourier (Grenoble) 12 (1962), 1--124.

\bibitem[7] I. Kol\'ar, P. Michor and J. Slov\'ak, {\it Natural Operations
in Differential Geometry}, Springer,
Berlin, 1993.

\bibitem[8] I.S. Krasilschik, A.M. Vinogradov and V.V. Lychagin, {\it
Geometry of Jet Spaces and Nonlinear 
Differential Equations}, Gordon and Breach, New York, 1986.

\bibitem[9] D. Krupka and J. Jany\u ska, {\it Lectures on Differential
Invariants}, Brno University, Czech Republic, 
1990.

\bibitem[10] D. Krupka, Local invariants of linear connections, Colloq.
Math. Soc. J\'anos Bolyai, 31. Differential 
Geometry, Budapest, 1979; North Holland, 1982, 349--369.

\bibitem[11] D. Krupka, Natural Lagrangian structures, Semester on Diff.
Geom., Banach Center, Warsaw, 
Sept.-Dec. 1979; Banach Center Publications 12, 1984, 185--210. 

\bibitem[12] M. Krupka, {\it Natural Operators on Vector Fields and Vector
Distributions}, Dr Dissertation, 
Faculty of Science, Masaryk University, Brno, Czech Republic, 1995

\bibitem[13] M. Krupka, Orientability of higher order Grassmannians, Math.
Slovaca 44 (1994), 107--115.

\bibitem[14] A.~Nijehuis, Natural bundles and their general properties,
Differential Geometry, in honour 
of K.~Yano, Kinokuniya, Tokyo, 1972, 317--334. 

\bibitem[15] J.M. Souriau, {\it Structure des Systemes Dynamiques}, Dunod,
Paris, 1970. 

\vskip 2cm

D.R.GRIGORE\footnote{Permanent Address: Department of Theoretical Physics,
Institute of Atomic Physics,
Bucharest-M\u agurele, P.O.Box MG6, Rom\^ania, e-mail:
grigore\@theor1.ifa.ro, grigore\@roifa.ifa.ro} 

D. KRUPKA

Department of Mathematics

Silesian University at Opava

Bezrucovo nam. 13, 74601 Opava

Czech Republic

e-mail: demeter.krupka\@fpf.slu.cz
\enddocument